\documentclass[english]{article}
\usepackage[T1]{fontenc}
\usepackage[latin1]{inputenc}
\usepackage{babel}
\usepackage{graphics}
\usepackage{psfrag}
 

\newtheorem{prop}{Proposition}[section]
\newtheorem{lemma}[prop]{Lemma}
\newtheorem{theorem}[prop]{Theorem}
\newtheorem{defi}[prop]{Definition}
\newtheorem{rem}[prop]{\em Remark\/}
\newtheorem{cor}[prop]{Corollary}
\newtheorem{cond}[prop]{Conditions}
\newtheorem{ex}[prop]{\em Example\/}

\newcommand{\bp}{\begin{prop}}
\newcommand{\bl}{\begin{lemma}}
\newcommand{\bt}{\begin{theorem}}
\newcommand{\bd}{\begin{defi}}
\newcommand{\br}{\begin{rem}}
\newcommand{\be}{\begin{equation}}
\newcommand{\bea}{\begin{eqnarray}}
\newcommand{\bpr}{\begin{proof}}
\newcommand{\bc}{\begin{cor}}
\newcommand{\bco}{\begin{cond}}
\newcommand{\bex}{\begin{ex}}

\newcommand{\ep}{\end{prop}}
\newcommand{\el}{\end{lemma}}
\newcommand{\et}{\end{theorem}}
\newcommand{\ed}{\end{defi}}
\newcommand{\er}{\end{rem}}
\newcommand{\ee}{\end{equation}}
\newcommand{\eea}{\end{eqnarray}}
\newcommand{\epr}{\end{proof}}
\newcommand{\ec}{\end{cor}}
\newcommand{\eco}{\end{cond}}
\newcommand{\eex}{\end{ex}}

\newcommand{\nn}{ \nonumber \\ }
\newcommand{\pr}{{\em Proof:\ }}
\newcommand{\qed}{\vrule height 5pt width 5pt depth 0pt}

\newcommand{\id}{{\rm id}}
 
\newcommand{\ui}{^{(1)}}
\newcommand{\uii}{^{(2)}}
\newcommand{\di}{_{(1)}}
\newcommand{\dii}{{}_{(2)}}
\newcommand{\ot}{\otimes}
\newcommand{\asr}{{A^R}}
\newcommand{\atr}{{^R\! A}}
\newcommand{\asl}{{_L\! A}}
\newcommand{\atl}{{A_L}}
\newcommand{\asrd}{{{\cal A}^*}}
\newcommand{\atrd}{{^*\!\! {\cal A}}}
\newcommand{\atld}{{{\cal A}_*}}
\newcommand{\asld}{{_*{\cal A}}}
\newcommand{\ci}{\circ}
\newcommand{\ila}{{\cal I}^L({\cal A})}

\renewcommand{\c}{{\cal C}}
\newcommand{\fsr}{{\ell_R}}
\newcommand{\ftr}{{_R\ell}}
\newcommand{\fsl}{{_L\ell}}
\newcommand{\ftl}{{\ell_L}}
\newcommand{\lu}{\leftharpoonup}
\newcommand{\ld}{\leftharpoondown}
\newcommand{\ru}{\rightharpoonup}
\newcommand{\rd}{\rightharpoondown}
\newcommand{\lsr}{\lambda\!^*}
\newcommand{\ltr}{^*\!\! \lambda}
\newcommand{\phisl}{_* \!\phi}
\newcommand{\phisr}{\phi\!^*}
\newcommand{\phitr}{^*\!\! \phi}
\newcommand{\phitl}{\phi_*}
\renewcommand{\i}{\iota}
\newcommand{\ib}{{\bar{\iota}}}
\newcommand{\x}{\times}

\newcommand{\stac}[1]{\stackrel{\ot}{_{_{#1}}}}

\newcommand{\err}{\Upsilon}

\newcommand{\setc}[1]{\setcounter{equation}{#1}}

\newcommand{\lb}{\label}
\begin{document}

\large
\title{\bf An alternative notion of Hopf Algebroid }
 
\author{\sc Gabriella B\"ohm $^1$ \\}

\date{}
 
\maketitle
\normalsize 
\footnotetext[1]{
Research Institute for Particle and Nuclear Physics, Budapest,
H-1525 Budapest 114, P.O.B. 49, Hungary\\
E-mail: BGABR@rmki.kfki.hu\\
  Supported by the Hungarian Scientific Research Fund, OTKA --
  T 034 512, FKFP -- 0043/2001 and the Bolyai J\'anos Fellowship}
\begin{abstract} In \cite{BKSz} a new notion of Hopf 
algebroid has been introduced. It was shown to be inequivalent to the
structure introduced under the same name in \cite{Lu}. We review this
new notion of Hopf algebroid. 
We prove that two Hopf algebroids are isomorphic as bialgebroids if
and only if 
their antipodes are related by a `twist' i.e. are deformed by the
analogue of a character.
A precise
relation to weak Hopf algebras is given. After the review of the
integral theory of Hopf algebroids we show how a right
bialgebroid can be made a Hopf algebroid in the presence of a
non-degenerate left integral. This can be interpreted as the `half of
the Larson-Sweedler theorem'. As an application we construct the Hopf
algebroid symmetry of an abstract depth 2 Frobenius extension \cite{FD2}.
\end{abstract}

\bigskip
\bigskip
\bigskip
\bigskip

%\Huge
%\[{Preliminary\ version}\]
\normalsize 
\vfill\eject

\section{Introduction}

Recently many authors introduced generalizations of bialgebras. These
structures are common in the feature that they do not need to be
algebras rather bimodules over some -- possibly non-commutative --
base ring $L$. In the paper \cite{BreMi} the notions of Lu's
bialgebroid \cite{Lu} (axiomatized in a more compact form in
\cite{Sz}), Xu's bialgebroid with anchor \cite{Xu,Xu1} and Takeuchi's
$\times_L$-bialgebra \cite{T} were shown to be equivalent. We use the
definition as follows:
\bd \lb{lbgd}
A {\em left bialgebroid} (or Takeuchi $\times_L$-bialgebra) ${\cal
A}_L$ consists of the data $(A,L,s_L,t_L,$\break$\gamma_L,\pi_L)$. The
$A$ and 
$L$ are associative unital rings, the total and base rings,
respectively. The $s_L:L\to A$ and $t_L:L^{op}\to 
A$ are ring homomorphisms such that the images of $L$ in $A$ commute
making $A$ an $L-L$ bimodule:
\be  l\cdot a\cdot l^{\prime}\colon = s_L(l) t_L(l^{\prime}) a .
\lb{elbim} \ee
The bimodule (\ref{elbim}) is denoted by ${_L A_L}$. The triple
$({_L A_L},\gamma_L,\pi_L)$ is a comonoid in  ${_L {\cal M}_L}$, the
category of $L-L$ bimodules. Introducing
the Sweedler's notation $\gamma_L(a)\equiv a\di \ot a\dii\in \atl \ot
\asl$ the identities
\bea  a\di t_L(l) \ot a\dii &=& a\di \ot a\dii s_L(l) \lb{cros}\\
      \gamma_L(1_A)&=& 1_A\ot 1_A \\
      \gamma_L(ab)&=&\gamma_L(a) \gamma_L(b) \lb{gmp} \\
      \pi_L(1_A) &=& 1_L \\
      \pi_L\left(a s_L\ci \pi_L(b)\right)=&\pi_L&(ab)=
       \pi_L\left(a t_L\ci \pi_L(b)\right)
\eea
are  required for all $l\in L$ and $a,b\in A$. The requirement 
(\ref{gmp}) makes sense in the view of (\ref{cros}).
\ed
With the help of the maps $s_L$ and $t_L$ we can introduce four
commuting actions of $L$ on $A$. They give rise to $L$-modules
\begin{eqnarray} &{_L A}:\ l\cdot a=s_L(l)a \qquad
     &{A_L}:\  a\cdot l=t_L(l)a \nonumber\\
     &{A^L}:\  a\cdot l=as_L(l) \qquad
     &{^L A}:\ l\cdot a=at_L(l). \nonumber
\end{eqnarray}
One defines the bimodules ${^L A^L}$, ${^L A_L}$ and ${_L A^L}$ in the
obvious way.

Throughout the paper it is a typical situation that the same ring $A$
carries different $L$-module structures. In this situation the usual
notation $A\stac{L}A$ is ambiguous. Our notation of bimodule tensor
products is explained at the beginning of Section 2.

In \cite{Lu} J. H. Lu introduced the notion of  {\em Hopf algebroid} 
as a triple $({\cal A}_L,S,\xi)$ consisting of a 
left bialgebroid ${\cal A}_L=(A,L,s_L,t_L,\gamma_L,\pi_L)$ such that
$A$ and $L$ are algebras over a commutative ring $k$. It is equipped
with an antipode $S:A\to A$. The $S$ is required to be an
anti-automorphism of the $k$-algebra $A$ satisfying 
\bea  S\ci t_L&=&s_L\lb{lu1}\\
m_A\ci (S\ot_L \id_A)\ci \gamma_L&=&t_L\ci \pi_L\ci S \lb{lu2}\\
m_A\ci (\id_A\ot_k S)\ci \xi\ci \gamma_L &=& s_L\ci \pi_L
\lb{lu3}\eea
where $m_A$  is the multiplication in $A$  and $\xi$ is  a
section of the canonical projection $p_L: A\ot_k A\to A_L \ot {_L A}$
that is $\xi$ is a map $A_L \ot{_L A}\to A\ot_k A$ 
satisfying $p_L\ci \xi = \id_{A\ot_L A}$.

\smallskip

Following the result of \cite{Szy,Lo} -- proving that an irreducible
finite index depth 2 
extension of von Neumann factors can be realized as a crossed product
with a finite dimensional $C^*$-Hopf algebra -- big effort has been
made in order to make connection between more and more general kinds
of extensions and of `quantum symmetries'
\cite{Oc,Da,EN,Ya,NW}. Allowing for {\em reducible} finite index D2
extensions of $II_1$ von 
Neumann factors in \cite{NV1} and of von Neumann algebras with finite
centers in \cite{NSzW} the symmetry of the extension was shown to be
described by a {\em finite dimensional $C^*$-weak Hopf algebra}
introduced in \cite{BSz,N,BNSz}. A Galois correspondence has been
established in \cite{NV2} in the case of finite index finite depth
extensions of $II_1$ factors. The infinite index D2 case has been
treated in \cite{EV,E}  for arbitrary von Neumann algebras endowed with
a regular operator valued weight.

In the paper
\cite{KSz} depth 2 extensions of {\em rings} have been
investigated. It was 
shown that there exists a canonical dual pair of (finite)  bialgebroids
associated to such a ring extension. 

Studying the bialgebroids corresponding to a depth 2 {\em
Frobenius} extension of rings one easily generalizes the
formula describing the antipode in \cite{Szy,Lo,NW} to an
anti-automorphism of the total ring satisfying the axioms
(\ref{lu1}-\ref{lu2}). This supports the expectation 
that in the case of a depth 2 Frobenius extension of rings the
canonical  
bialgebroids obtained in \cite{KSz} can be made Hopf
algebroids. However, no effort, made for checking the Lu-axioms
\cite{Lu} in this situation, 
brought success. (See however \cite{K1,K2} studying interesting
subcases.) As a matter of fact the section $\xi$, appearing in  
the definition in \cite{Lu}, does not come naturally into this
context. 

This lead us to the introduction of an alternative notion of Hopf
algebroid in \cite{BKSz} -- the prototype of which is the symmetry of
a depth 2 Frobenius extension of rings.
In \cite {BKSz} three equivalent sets of 
axioms are formulated.\footnote{added in proof: The final version of
\cite{BKSz} contains four 
equivalent definitions. A 'zeroth' one has been added later.}
The first definition is
analogous to the one in \cite{Lu} in the aspect that a left
bialgebroid ${\cal A}_L= 
(A,L,s_L,t_L,\gamma_L,\pi_L)$ is equipped
with a bijective antipode map $S$ which is a ring anti-automorphism and
relates the left and right
$L$-module structures as in (\ref{lu1}). Also the antipode axiom
(\ref{lu2}) 
is identical. The main difference is that in the definition of
\cite{BKSz} no reference to a section `$\xi$' is needed. In its stead
we deal with the maps  
\be (S\ot S)\ci \gamma_L^{op}\ci S^{-1}\qquad {\rm and} \qquad
(S^{-1}\ot S^{-1})\ci \gamma_L^{op}\ci S\lb{gr}.\ee 
In the Hopf algebra case both are equal to the
coproduct itself. In this more general case 
the images of $A$ under $\gamma_L$ and the maps (\ref{gr}) are
however different. We require the two maps of (\ref{gr}) to be 
equal and be both a left and a right comodule map.
This form of the axioms does not contain a second antipode axiom.
This definition is cited at the end of Section \ref{def} of this
paper. 

Analyzing the consequences of this first  definition  one observes
a hidden right bialgebroid
\footnote{ For the definition of the right bialgebroid \cite{KSz} see
Definition \ref{rbgd} below.}  
structure on $A$ with the coproduct given
by the equal maps (\ref{gr}). This lead us to the second  `symmetric'
definition in \cite{BKSz} where the two antipode axioms have analogous
forms using both the left and right bialgebroid structures of $A$ in a
symmetric way. We cite this definition in Definition \ref{haid} below.
 
The third definition in \cite{BKSz}
is formulated without the explicit use of the antipode map. It
borrows the philosophy of \cite{Sch} where   
a Hopf algebroid like object (possibly without antipode), the so called
$\times_L${\em -Hopf algebra} was introduced as  a left bialgebroid
s.t. the map 
\[ \alpha: {^L A}\stackrel{\ot}{_{ _{L^{op}}}} A_L\to A_L\ot {_L A}\qquad 
a\ot b \mapsto a\di \ot a\dii b\]
is bijective. 
Requiring however the bijectivity of the maps $\alpha$ and
its co-opposite $\beta$ (leading to a bijective antipode in the 
bialgebra case) for a left bialgebroid does not imply the existence of  
an antipode. In order to 
have a definition which is equivalent to the other two both related
bialgebroid structures are needed. This definition explicitly shows
that the Hopf algebroids in the sense of \cite{BKSz} are
$\times_L$-Hopf algebras. 
  
It is proven in \cite{BKSz} that (weak) Hopf algebras with bijective
antipode are Hopf
algebroids in the sense of \cite{BKSz}. Also some examples of Lu-Hopf
algebroids \cite{Lu, BreMi, KhalRang, ConnMo2} are shown to be
examples.  

The two notions of Hopf algebroid -- the one introduced in \cite{Lu}
and the one in \cite{BKSz} --  were shown to be inequivalent by giving
an example of the Hopf algebroid \cite{BKSz} that does not satisfy the
axioms of \cite{Lu}. This example is discussed in the Section
\ref{ex} of this paper in more detail.

On the other hand no example of Lu-Hopf algebroid is known to us at
the moment that does not satisfy the axioms of \cite{BKSz}. This
leaves open the logical possibility that the Lu-Hopf algebroid was a
subcase of the one introduced in \cite{BKSz}. Until now we could
neither prove nor exclude by examples this possibility.

In \cite{BKSz} the theory of non-degenerate integrals in a Hopf
algebroid is developed. Though the axioms of the Hopf algebroid in
\cite{BKSz} are by no means self-dual, it is proven in \cite{BKSz}
that if there exists a non-degenerate left integral $\ell$ in a Hopf
algebroid ${\cal A}$ then its dual (with respect to the base ring)
also carries a Hopf algebroid structure which is unique upto an
isomorphism of bialgebroids. The dual of the bialgebroid isomorphism
class of ${\cal A}$ is defined as the bialgebroid isomorphism
class of the Hopf algebroid constructed on the dual ring.

In this paper we present a review and also some new results on Hopf
algebroids. In the 
Section \ref{bgd} we review some results about bialgebroids that were
obtained in the papers \cite{Sw,T,Sz,KSz}. In the Section \ref{def} we
repeat the definition of the Hopf algebroid given in \cite{BKSz} and
cite some basic results without proofs. In the Section \ref{ex} we
generalize the notion of the {\em twist 
of a Hopf algebra} -- introduced in \cite{ConnMo} -- to Hopf
algebroids. In particular we prove that twisted (weak) Hopf algebras
are Hopf algebroids in the sense of \cite{BKSz}. A most important
example is the Connes-Moscovici algebra ${\cal 
H}_{FM}$ \cite{ConnMo2}. By twisting cocommutative Hopf algebras we construct
examples that do not satisfy the Hopf algebroid axioms of
\cite{Lu}. We give a sufficient and necessary criterion on a Hopf
algebroid under which it is  
a (twisted version) of weak Hopf algebra. In the Section \ref{int} we
review the 
integral theory of Hopf algebroids from \cite{BKSz} without
proofs. The Section \ref{LS} deals with the question how can we make
a right bialgebroid with a non-degenerate left integral into a Hopf
algebroid. The result of Section \ref{LS} can be interpreted as the
generalization of the `easier half of the Larson-Sweedler theorem'
\cite{LaSw}. 

\section{Bialgebroids}
\lb{bgd}
\setc{0}
The total ring of a bialgebroid carries different module structures
over the base ring. In this situation we make the following notational
convention. In writing module tensor products we write out the (bi-)
module factors explicitly. For the $L$-module tensor product of the
bimodules ${_L A^L}$ and ${_L A_L}$, for example we write ${_L
A^{L^{\prime}}}\stac{L^{\prime}}{_{L^{\prime}} A_L}$, where $L^{\prime}$
stands for another copy of $L$, and it has been introduced to show
explicitly which module structures are involved in the tensor product.

In order not to make the formulas more complicated than necessary we
make a further simplification. In those situations in which it is
clear from the tensor factors themselves over which ring the tensor
product is taken, we do not denote it under the symbol
$\otimes$. I.e. for the $L$-module tensor product of the right
$L$-module $A_L$ and the left $L$-module ${_L A}$, for example, we
write $A_L\otimes _LA$.

The {\em bialgebroid} \cite{Lu,Xu,Sz} or -- what is equivalent to it -- a
Takeuchi $\x_L$-bialgebra \cite{T}  is a generalization of the
bialgebra in the sense that it is not an algebra rather a 
bimodule over a non-commutative ring $L$. We use Definition \ref{lbgd}
of the left bialgebroid.
We use the name {\em left} bialgebroid as in \cite{KSz} since the
`opposite structure' was called a {\em right bialgebroid} in
\cite{KSz}: 
\bd \lb{rbgd}
A {\em right bialgebroid} ${\cal A}_R$ consists of the data
$(A,R,s_R,t_R,\gamma_R,\pi_R)$. The $A$ and 
$R$ are associative unital rings, the total and base rings,
respectively. The $s_R:R\to A$ and $t_R:R^{op}\to 
A$ are ring homomorphisms such that the images of $R$ in $A$ commute
making $A$ an $R-R$ bimodule:
\be  r\cdot a\cdot r^{\prime}\colon = a s_R(r^{\prime}) t_R(r).
\lb{erbim}\ee
The bimodule (\ref{erbim}) is denoted by ${^R\! A^R}$. The triple
$({^R\! A^R},\gamma_R,\pi_R)$ is a comonoid in ${_R {\cal
M}_R}$. Introducing 
the Sweedler's notation $\gamma_R(a)\equiv a\ui \ot a\uii\in \asr \ot
\atr$ the identities
\bea  s_R(r)a\ui  \ot a\uii &=& a\ui \ot  t_R(r) a\uii \nn
      \gamma_R(1_A)&=& 1_A\ot 1_A \nn
      \gamma_R(ab)&=&\gamma_R(a) \gamma_R(b) \nn
      \pi_R(1_A) &=& 1_R \nn
      \pi_R\left(s_R\ci \pi_R(a) b \right)=&\pi_R&(ab)=
       \pi_R\left(t_R\ci \pi_R(a) b \right)
\nonumber\eea
are  required for all $r\in R$ and $a,b\in A$.
\ed
In addition to the bimodule ${^RA^R}$ we introduce also
\[  {_RA_R}:\qquad r\cdot a\cdot r^{\prime}\colon = s_R(r) t_R(
r^{\prime}) a. \]
If ${\cal A}_L=(A,L,s_L,t_L,\gamma_L,\pi_L)$ is a left bialgebroid
then so is the co-opposite 
$({\cal A}_L)_{cop}=(A,L^{op},t_L,s_L,\gamma_L^{op},$ $\pi_L)$ -- where
$\gamma_L^{op}:A\to \asl\ot \atl$ maps $a$ to $a\dii\ot a\di$ --. The
opposite  
$({\cal A}_L)^{op}=(A^{op},L,t_L,s_L,\gamma_L,\pi_L)$ is a right
bialgebroid. 

We use the terminology of homomorphisms of bialgebroids as introduced
in \cite{Sz2}:
\bd A {\em left bialgebroid homomorphism} ${\cal A}_L \to {\cal
A}^{\prime}_{L^{\prime}}$ is a pair of ring homomorphisms $(\Phi:A\to
A^{\prime}, 
\phi: L\to L^{\prime})$ such that 
\bea s^{\prime}_L\ci \phi &=& \Phi\ci s_L\nn
     t^{\prime}_L\ci \phi &=& \Phi\ci t_L\nn
     \pi^{\prime}_L\ci \Phi &=& \phi\ci \pi_L\nn
     \gamma_L^{\prime}\ci \Phi&=& (\Phi\ot \Phi)\ci \gamma_L.
\nonumber\eea
The last condition makes sense since by the first two
conditions $\Phi\ot \Phi$ is a well defined map $\atl \ot \asl\to
{A^{\prime}_{L^{\prime}}}\ot {_{L^{\prime}}}A^{\prime}$.
The pair $(\Phi,\phi)$ is an {\em isomorphism of left bialgebroids} if it is
a left bialgebroid homomorphism such that both $\Phi$ and $\phi$ are 
bijective. 

A right bialgebroid homomorphism (isomorphism) ${\cal A}_R\to {\cal
A}^{\prime}_{R^{\prime}}$ is a left bialgebroid homomorphism
(isomorphism) $({\cal 
A}_R)^{op}\to ({\cal A}^{\prime}_{R^{\prime}})^{op}$.
\ed

\smallskip

Let ${\cal A}_L$ be a left bialgebroid. The equation (\ref{elbim})
describes two $L$-modules $\atl$ and $\asl$. Their $L$-duals
are the additive groups of $L$-module maps:
\[   \atld\colon = \{ \phitl: \atl \to {L_ L} \} \quad {\rm and} \quad
     \asld\colon = \{ \phisl: \asl \to {_L L} \} \]
where ${_L L}$ stands for the left regular and $L_L$ for the right
regular $L$-module. Both $\atld$ and $\asld$ carry left $A$ module
structures via the transpose of the right regular action of $A$. For
$\phitl\in \atld, \phisl\in\asld$ and $a,b\in A$ we have:
\[   \left(a\ru \phitl\right)(b) =\phitl(ba)\quad {\rm and}\quad
     \left(a\rd \phisl\right)(b) =\phisl(ba). \]
Similarly, in the case of a right bialgebroid ${\cal A}_R$ --  denoting
the left and right regular $R$-modules by $^R R$ 
and $R^R$, respectively, -- the two $R$-dual additive groups
\[  \asrd\colon = \{ \phisr: \asr \to {R^R} \} \quad {\rm and}\quad
    \atrd\colon = \{ \phitr: \atr \to {^R R} \} \]
carry right $A$-module structures:
\[ \left(\phisr \lu a \right)(b) =\phisr(ab) \quad {\rm and}\quad
   \left(\phitr \ld a \right)(b) =\phitr(ab). \]
The comonoid structures can be transposed to give monoid (i.e. ring)
structures to the duals. In the case of a left bialgebroid ${\cal A}_L$ 
\be  \left(\phitl {\psi_*} \right)(a)=
{\psi_*}\left( s_L \ci \phitl (a\di) a\dii \right)\quad {\rm and}\quad
\left(\phisl {_* \psi} \right)(a)=
{_* \psi}\left( t_L \ci \phisl (a\dii) a\di \right) \lb{ldual}\ee
for $\phisl, {_*\psi}\in \asld$, $\phitl,{\psi_*}\in \atld$ and $a\in
A$. 
Similarly, in the case of a right bialgebroid ${\cal A}_R$ 
\be \left(\phisr {\psi^*} \right)(a)=
\phisr\left( a\uii t_R \ci {\psi^*} (a\ui)\right)\quad {\rm and}\quad
 \left(\phitr {^*\! \psi} \right)(a)=
\phitr \left(a\ui  s_R \ci {^* \!\psi} (a\uii)\right)\lb{rdual}\ee
for $\phisr, {\psi^*}\in \asrd$, $\phitr,{^*\psi}\in \atrd$ and $a\in
A$.
In the case of a left bialgebroid ${\cal A}_L$  also the ring $A$ has 
right $\atld$- and   right $\asld$- module structures:
\[ a\lu \phitl = s_L\ci \phitl (a\di)a\dii \quad {\rm and}\quad
     a\ld \phisl = t_L\ci \phisl (a\dii)a\di \]
for $\phitl\in \atld$, $\phisl\in \asld$ and $a\in A$.

Similarly, in the case of a right bialgebroid ${\cal A}_R$  the ring
$A$ has left $\asrd$- and left $\atrd$ structures:
\[ \phisr\ru a = a\uii t_R\ci \phisr(a\ui )\quad {\rm and}\quad
     \phitr\rd a = a\ui  s_R\ci \phitr(a\uii) \]
for $\phisr\in \asrd$,  $\phitr\in \atrd$ and $a\in A$.

It is also proven in \cite{KSz} that if the $L$ ($R$) module structure
on $A$ is finitely generated projective then the corresponding dual
has also a bialgebroid structure.

\section{Hopf Algebroid}
\lb{def}
\setc{0}

Let ${\cal A}_L=(A,L,s_L,t_L,\gamma_L,\pi_L)$ be a left bialgebroid
and  ${\cal A}_R=(A,R,s_R,t_R,\gamma_R,\pi_R)$ a right bialgebroid
such that the base rings are anti-isomorphic $R\simeq L^{op}$. Require
that
\be  s_L(L)=t_R(R)\qquad {\rm and} \qquad t_L(L)=s_R(R) \lb{defi}\ee
as subrings of $A$. The requirement (\ref{defi}) implies that the
coproduct $\gamma_L$ is a quadro-module map 
${_L^R\! A^R_L}\to {_L^R A
_{L^{\prime}}}\stackrel{\ot}{_{_{L^{\prime}}}} 
{_{L^{\prime}}A^R_L}$
and $\gamma_R$ is a quadro-module map 
${_L^R A^R _L}\to 
{_L^R A ^{R^{\prime}}}\stackrel{\ot}{_{_{R^{\prime}}}}
{^{R^{\prime}}A^R_L}$. (The $L^{\prime}$ and $R^{\prime}$ denote
another copy of $L$ and $R$, respectively). This allows us to require that 
\bea  &(\gamma_L \ot \id_A)\ci \gamma_R &= (\id_A \ot \gamma_R)\ci
\gamma_L \nn
      &(\gamma_R \ot \id_A)\ci \gamma_L &= (\id_A \ot \gamma_L)\ci
\gamma_R
\lb{defii}
\eea
as maps $A\to \atl \stackrel{\ot}{_{_L}} {\asl^R}
\stackrel{\ot}{_{_R}} {\atr} $ and
$A\to \asr \stackrel{\ot}{_{_R}} {\atr_L}
\stackrel{\ot}{_{_L}} \asl$,
respectively. 

Let  $S:A\to A$ be a bijection of additive groups such that 
$S$ is a twisted  isomorphism of bimodules ${^L A_L}\to {_L A^L}$ and ${^R
A_R}\to {_R A ^R}$ that is
\be    S(t_L(l)a t_L(l^{\prime}))=s_L(l^{\prime})S(a) s_L(l) \qquad 
            S(t_R(r^{\prime})a t_R(r))=s_R(r)S(a) s_R(r^{\prime}) 
\lb{defiii} \ee 
for all $l,l^{\prime}\in L$, $r,r^{\prime}\in R$ and $a\in A$.
The requirement
(\ref{defiii}) makes the expressions $S(a\di)a\dii$ and $a\ui S(a\uii)$
meaningful. We require
\be      S(a\di)a\dii =s_R\ci \pi_R (a) \qquad 
         a\ui S(a\uii)=s_L\ci \pi_L (a) \lb{defiv}
\ee
for all $a$ in $A$.

\bd \lb{haid}
The triple ${\cal A}=({\cal A}_L,{\cal A}_R,S)$ satisfying
(\ref{defi}), (\ref{defii}), (\ref{defiii}) and (\ref{defiv}) is a
{\em Hopf algebroid}.
\ed 

Throughout the paper we use the analogue of the Sweedler's notation:
$\gamma_L(a)=a\di \ot a\dii$ and $\gamma_R(a)=a\ui\ot a\uii$.

For a Hopf algebroid ${\cal A}$ also the opposite  ${\cal
A}^{op}=({\cal A}_R^{op},{\cal A}_L^{op}, S^{-1})$ and the co-opposite
${\cal A}_{cop}=({\cal A}_{L\ cop}, {\cal A}_{R\ cop},S^{-1})$ are
Hopf algebroids.

In the Section \ref{ex} we will investigate the way in which (weak)
Hopf algebra is a  subcase.

The antipode of a Hopf algebra is a bialgebra anti-homomorphism. This
property generalizes to Hopf algebroids as  
\bp \lb{Sisom} Both
$(S,\pi_R\ci s_L)$ and $(S^{-1},\pi_R\ci t_L)$ are left bialgebroid
isomorphisms ${\cal A}_L\to ({\cal A}_R)^{op}_{cop}$. That is 
\bea &s_R\ci \pi_R\ci s_L=S\ci s_L \qquad
     &s_R\ci \pi_R\ci t_L=S^{-1}\ci s_L \nn
     &t_R\ci \pi_R\ci s_L=S\ci t_L \qquad
     &t_R\ci \pi_R\ci t_L=S^{-1}\ci t_L \nn
     &\pi_R\ci s_L\ci \pi_L=\pi_R\ci S\qquad
     &\pi_R\ci t_L\ci \pi_L=\pi_R\ci S^{-1}\nn
     &S_{A\ot_L A}\ci \gamma_L=\gamma_R\ci S\qquad
     &S^{-1}_{A\ot_R A}\ci \gamma_L=\gamma_R\ci S^{-1}
\nonumber\eea
where $S_{A\ot_LA}$ is a map $\atl \ot \asl \to \asr\ot \atr$, it maps
$a\ot b$ to $S(b)\ot S(a)$. Similarly,  $S_{A\ot_R A}$ is a map
$\asr\ot \atr \to \atl \ot \asl$, it maps  $a\ot b$ to $S(b)\ot S(a)$.
\ep

The datum $({\cal A}_L,{\cal A}_R,S)$ determining a Hopf algebroid is
somewhat redundant. Indeed, suppose that we have  given only a left
bialgebroid ${\cal A}_L=(A,L,s_L,t_L,\gamma_L,\pi_L)$ and an
anti-isomorphism $S$ of the 
total ring $A$ satisfying
\bea S\ci t_L &=& s_L\lb{lui}\\
%S(a\di) a\dii &=& t_L\ci \pi_L \ci S(a) \lb{luii}\\
m_A\ci(S\ot \id_A)\ci \gamma_L&=&t_L\ci\pi_L\ci S\\
%\gamma_L\ci S^2&=& (S^2\ot S^2)\ci \gamma_L \lb{luiii}\\
S_{A\ot_L A}\ci \gamma_L\ci S^{-1} &=&
S^{-1}_{A\ot_R A}\ci \gamma_L\ci S\lb{luiii}\\
(\gamma_L\ot \id_A)\ci \gamma_R = (\id_A \ot \gamma_R)\ci
\gamma_L \quad &\ &\quad
(\gamma_R\ot \id_A)\ci \gamma_L = (\id_A \ot \gamma_L)\ci
\gamma_R \lb{luiv}
\eea
for $m_A$ the multiplication in $A$ and  $\gamma_R\colon = S_{A\ot_L A}\ci
\gamma_L\ci S^{-1}$. Then the right 
bialgebroid ${\cal A}_R$ -- together with which $({\cal A}_L,{\cal
A}_R,S)$ is a Hopf algebroid -- can be reconstructed upto a trivial
bialgebroid isomorphism. Namely, it follows from Proposition
\ref{Sisom} that 
\[{\cal A}_R= (A,R,S\ci s_L\ci \nu^{-1},  s_L\ci \nu^{-1}, S_{A\ot_L
A}\ci\gamma_L\ci S^{-1},\nu \ci \pi_L\ci S^{-1})\]
for an arbitrary  isomorphism $\nu:L^{op}\to R$.

\section{Twist of the Hopf algebroid}
\lb{ex}
\setc{0}

It is clear that being given the left and right bialgebroids ${\cal
A}_L$ and ${\cal A}_R$ satisfying the axioms (\ref{defi}) and
(\ref{defii}) the antipode -- if exists -- is unique. Indeed, if both
$S_1$ and $S_2$ make $({\cal A}_L,{\cal A}_R,S_1)$ and $({\cal
A}_L,{\cal A}_R,S_2)$ Hopf algebroids then
\bea S_2(a)&=& s_R\ci\pi_R(a\ui) S_2(a\uii)= 
S_1({a\ui}\di) {a\ui}\dii S_2 (a\uii)=
S_1({a\di}) {a\dii}\ui S_2 ({a\dii}\uii)=\nn
&=&S_1(a\di) s_L\ci\pi_L (a\dii) = S_1(a).
\nonumber\eea
There are some examples however in which only
the left bialgebroid structure is naturally given and we have some
ambiguity in the choice of the right bialgebroid structure and the
corresponding antipode. (See for example the Hopf algebroid symmetry
of a depth 2 Frobenius extension of rings in Section 3 of \cite{BKSz}
or at the end of Section \ref{LS} below. In this example the
ambiguity is 
nicely controlled by the Radon-Nykodim derivative relating the possible
Frobenius maps.) In the following we address the question more
generally: Given a left bialgebroid ${\cal A}_L$ how are the
possible antipodes satisfying the conditions (\ref{lui}-\ref{luiv}) related? 
\bd Let $({\cal A}_L,S)$ be subject to the conditions
(\ref{lui}-\ref{luiv}). An 
invertible element $g_*$ of $\atld$ is called a {\em twist} of $({\cal
A}_L,S)$ if for all elements $a,b$ of $A$
\bea
&i)& 1_A\lu g_* = 1_A \nn
&ii)& (a\lu g_*)(b \lu g_*) = ab\lu g_* \nn
&iii)&  S(a\di)\lu g_* \ot a\dii = S(a\di) \ot a\dii\lu g_*^{-1}.  
\lb{twistdef}
\eea
The condition iii) is understood to be an identity in the product of
the modules $A^L$ and the left $L$-module on $A$:
\[ l\cdot a = s_L\ci g_*^{-1} \ci s_L(l) a \]
where $g_*\ci s_L$ is an automorphism of $L$ with inverse $g_*^{-1}\ci
s_L$. 
\ed
The {\em twist} of Definition \ref{twistdef} generalizes the notion of
the character on a Hopf algebra.
Clearly the twists of $({\cal A}_L,S)$ form a group.
\bt \lb{twist}
Let $({\cal A}_L,S)$ be subject to the conditions (\ref{lui}-\ref{luiv}). Then
$({\cal A}_L,S^{\prime})$ is subject to the conditions
(\ref{lui}-\ref{luiv})  if and only
if there exists a twist $g_*$ of $({\cal A}_L,S)$ such that $S^{\prime}(a)=
S(a\lu g_*)$ for all $a\in A$.
\et

\br For $({\cal A}_L,S)$ a {\em Hopf algebra} the twisted 
antipode of the above form was introduced in \cite{ConnMo}.
In the view of Theorem \ref{twist} the twisted Hopf algebras in
\cite{ConnMo} are Hopf algebroids in the sense of Definition
\ref{haid}. 
\er

{\em Proof of the Theorem: } {\em if part:}
For a twist $g_*$ the map $S_g(a) \colon = S( a\lu g_*)$ is bijective
with inverse $S_g^{-1}(a) = S^{-1}(a) \lu g_*^{-1}$. It is
anti-multiplicative by ii) of (\ref{twistdef}).
Using the property i) of (\ref{twistdef}) 
\[ S_g\ci t_L(l)= S\left( t_L(l)\lu g_* \right) = S\left((1_A\lu g_*)t_L(l)\right)=
S\ci t_L(l)=s_L(l). \]
By $\gamma_L(a\lu g_*)=(a\di \lu g_*) \ot a\dii$ and $S\ci t_R\ci
\pi_R=t_L\ci \pi_L\ci S$ we have
\[ S_g(a\di) a\dii = S(a\di \lu g_*) a\dii= s_R\ci \pi_R(a\lu g_*)=
t_L\ci \pi_L \ci S_g(a).\]
In order to check the property (\ref{luiii}) of $({\cal A}_L,S_g)$
rewrite iii) of (\ref{twistdef}) into the equivalent form 
\be S(a\di)\ot S_g(a\dii) = S_g^{-1}\ci S^2(a\di)\ot S(a\dii).
\lb{iiipr}\ee
Then using the fact that $({\cal A}_L,S)$ satisfies (\ref{luiii})
introduce $a\ui\ot a\uii \colon = S_{A\ot_L A}\ci
\gamma_L\ci S^{-1}(a)\equiv S^{-1}_{A\ot_R A} \ci \gamma_L\ci
S(a)$. By (\ref{iiipr}) we have 
\bea &S_g& (S_g^{-1}(a)\dii)\ot  S_g(S_g^{-1}(a)\di)=
 S_g(S^{-1}(a)\dii) \ot  S_g(S^{-1}(a)\di\lu g_*^{-1})=\nn
&=&  S_g(S^{-1}(a)\dii) \ot  S(S^{-1}(a)\di)=
S(S^{-1}(a)\dii) \ot  S_g^{-1}\ci S^2(S^{-1}(a)\di)=\nn
%&=& S^{-1}(S(a)\dii)\ot  S_g^{-1}(S(a)\di)=\nn
%&=&S_g^{-1}\left[S\left(S^{-1}[S(a)\dii]\lu g_*\right)\right]\ot  
%S_g^{-1}(S(a)\di)=\nn
%&=&S_g^{-1}\ci S(a\ui\lu g_*)\ot S_g^{-1}\ci S(a\uii)=\nn
%&=&S_g^{-1}\ci S((a\lu g_*)\ui)\ot S_g^{-1}\ci S((a\lu g_*)\uii)=\nn
%%%%%%%%%%%%%%%%%%%%%%%%%%%%%%%%%%%%%%%%%%%%%%%%%%%%%%%%%%%%%%%%%%%
&=&a\ui \ot S_g^{-1}\ci S(a\uii)=
 S_g^{-1}\ci S (a\ui \lu g_*)\ot  S_g^{-1}\ci S (a\uii)=\nn
%%%%%%%%%%%%%%%%%%%%%%%%%%%%%%%%%%%%%%%%%%%%%%%%%%%%%%%%%%%%%%%%%%%
&=&S_g^{-1}\left(S(a\lu g_*)\dii\right)\ot S_g^{-1}\left(S(a\lu
g_*)\di\right)= 
S_g^{-1}\left(S_g(a)\dii\right)\ot S_g^{-1}\left(S_g(a)\di\right).
\nonumber\eea
The last condition (\ref{luiv}) on $({\cal A}_L,S_g)$ follows then
easily by using the  two forms of
$S_{g \ A\ot_L A}\ci \gamma_L\ci S_g^{-1}(a)\equiv \gamma_{gR}(a)=
a\ui\ot S_g^{-1}\ci S(a\uii)$ 
and $\gamma_{gR}(a)= S_g\ci S^{-1}(a\ui)\ot a\uii$, respectively:
\bea  (\gamma_L\ot {\id_A})\ci \gamma_{gR}(a)&=&
{a\ui}\di\ot {a\ui}\dii\ot S_g^{-1}\ci S(a\uii)=
a\di \ot {a\dii}\ui \ot   S_g^{-1}\ci S({a\dii}\uii)=\nn
&=&(\id_A\ot \gamma_{gR})\ci \gamma_L(a)\nn
({\id_A}\ot \gamma_L)\ci \gamma_{gR}(a)&=& 
S_g\ci S^{-1}(a\ui)\ot {a\uii}\di \ot {a\uii}\dii=
S_g\ci S^{-1}({a\di}\ui)\ot {a\di}\uii\ot a\dii=\nn
&=&(\gamma_{gR} \ot \id_A)\ci \gamma_L(a).
\nonumber\eea
{\em only if part:}
Let  $({\cal A}_L,S)$ and $({\cal A}_L,S^{\prime})$ be subject to the
conditions (\ref{lui}-\ref{luiv}). By the considerations at the end of
Section \ref{def} we can construct the right
bialgebroids ${\cal A}_{R}$ and ${\cal A}_{R^{\prime}}^{\prime}$ such
that both 
${\cal A}=({\cal A}_L,{\cal A}_{R},S)$ and ${\cal A}^{\prime}=({\cal
A}_L,{\cal A}_{R^{\prime}}^{\prime},S^{\prime})$ are Hopf
algebroids. Denoting 
the coproducts in ${\cal A}_R$ and  ${\cal A}_{R^{\prime}}^{\prime}$ by
$\gamma_R(a)=a\ui\ot a\uii$ and
$\gamma_{R^{\prime}}^{\prime}(a)=a^{(1)^{\prime}}\ot a^{(2)^{\prime}}$,
respectively, the  maps
\bea \alpha:{^L A}\ot {A_L}\to \atl\ot \asl &\qquad & a\ot b\mapsto
a\di \ot a\dii b\nn
\beta: {A^L}\ot {_L A}\to {_L A}\ot \atl  &\qquad & a\ot b\mapsto a\dii
\ot a\di b 
\nonumber \eea
are easily shown to be bijections with inverses
\bea a\ui \ot S(a\uii)b = &\alpha^{-1}(a\ot b) &= a^{(1)^{\prime}}
\ot S^{\prime}(a^{(2)^{\prime}})b\nn 
     a\uii \ot S^{-1} (a\ui)b =& \beta^{-1} (a\ot b) &=
a^{(2)^{\prime}} \ot S^{\prime -1}(a^{(1)^{\prime}})b.
\lb{schinvun}
\eea
(This means that ${\cal A}_L$ and $({\cal A}_L)_{cop}$ are $\x_L$-Hopf
algebras in the sense of \cite{Sch}. )

We construct the twist $\pi_L\ci S^{-1} \ci S^{\prime}$. It is an
invertible element of $\atld$ with inverse $\pi_L\ci S^{\prime -1} \ci S$:
\bea [(\pi_L \ci  S^{-1} \ci S^{\prime})&&
\!\!\!\!\!\!\!\!\!\! \!\!\!\!\!
(\pi_L\ci S^{\prime -1} \ci S)](a) =
\pi_L \ci  S^{\prime -1} \ci S \left(
s_L\ci \pi_L \ci  S^{-1} \ci S^{\prime}(a\di) a\dii \right) =\nn
&=& \pi_L \ci  S^{\prime -1} \ci S \left(
s_L\ci \pi_L \ci  S^{-1}[S^{\prime}(a)^{(2)^{\prime}}] 
S^{\prime -1}[S^{\prime}(a)^{(1)^{\prime}}] \right) =\nn
&=& \pi_L \ci  S^{\prime -1} \ci S \left(
s_L\ci \pi_L \ci  S^{-1}[S^{\prime}(a)^{(2)}] 
S^{-1}[S^{\prime}(a)^{(1)}] \right) =\nn
&=& \pi_L \ci  S^{\prime -1} \left( S^{\prime}(a)^{(1)} 
s_R\ci\pi_R (S^{\prime} (a)\uii) \right) =\pi_L(a)
\lb{twistinv}
\eea 
where in the third step (\ref{schinvun}) has been used. The relation 
$(\pi_L\ci S^{\prime -1} \ci S)(\pi_L \ci  S^{-1} \ci S^{\prime})=\pi_L$
follows by interchanging the roles of $S$ and $S^{\prime}$.

For all elements $a,b$ in $A$ we have 
\bea &&1_A\lu \pi_L\ci S^{-1}\ci S^{\prime} = 
s_L\ci \pi_L\ci S^{-1}\ci S^{\prime} (1_A)=1_A\nn
&&(a\lu \pi_L\ci S^{-1}\ci S^{\prime})
(b\lu\pi_L\ci S^{-1}\ci S^{\prime })=\nn 
&&\qquad =s_L\ci \pi_L\ci S^{-1}\ci S^{\prime }(a\di)a\dii 
s_L \ci \pi_L\ci S^{-1}\ci S^{\prime }(b\di) b\dii=\nn
%&&\qquad =s_L\ci \pi_L\ci S^{-1}\ci S^{\prime} 
%\left(a\di t_L\ci \pi_L\ci
%S^{-1}\ci S^{\prime}(b\di)\right) a\dii b\dii=\nn
&&\qquad =s_L\ci \pi_L\left(
S^{-1}\ci S^{\prime}(a\di) 
t_L\ci \pi_L\ci S^{-1}\ci S^{\prime}(b\di)\right) a\dii b\dii=\nn
&&\qquad =s_L\ci \pi_L\left(S^{-1}\ci S^{\prime}(a\di)  
S^{-1}\ci S^{\prime}(b\di)\right) a\dii b\dii= 
s_L\ci  \pi_L\ci S^{-1}\ci S^{\prime}(a\di b\di)a\dii b\dii=\nn
&&\qquad = ab\lu  \pi_L\ci S^{-1}\ci S^{\prime}.
\nonumber\eea
Using (\ref{schinvun}) we can show that
\bea S(a&\lu& \pi_L \ci S^{-1}\ci S^{\prime})=
%S\left( s_L\ci \pi_L \ci S^{-1}\ci S^{\prime}(a\di)a\dii\right)=\nn&=&
S\left( s_L\ci \pi_L \ci S^{-1}[S^{\prime}(a)^{(2)^{\prime}}]
S^{\prime -1}[S^{\prime}(a)^{(1)^{\prime}}]\right)=\nn
&=&S\left( s_L\ci \pi_L \ci S^{-1}[S^{\prime}(a)^{(2)}]
S^{-1}[S^{\prime}(a)^{(1)}]\right)=
%\nn&=&
S^{\prime}(a)^{(1)} s_R\ci\pi_R (S^{\prime}(a)^{(2)})= 
S^{\prime}(a).
\lb{sg}
\eea
In order to check that $ \pi_L\ci S^{-1}\ci S^{\prime}$  satisfies iii)
of (\ref{twistdef}) rewrite (\ref{schinvun}) into the equivalent forms:
\[\begin{array}{crcl}
&a\ui \ot S^{\prime -1}\ci S(a\uii) = a^{(1)^{\prime}} &\ot&
a^{(2)^{\prime}} = S^{\prime}\ci S^{-1} (a\ui) \ot a\uii \nn
\Leftrightarrow& S^{-1}(S(a)\dii)\ot  S^{\prime -1}(S(a)\di)&=&
S^{\prime}(S^{-1}(a)\dii)\ot S(S^{-1}(a)\di)\nn
\Leftrightarrow&S(a\dii)\ot S^{\prime -1}\ci S^{2}(a\di)&=& 
S^{\prime}(a\dii)\ot S(a\di).\end{array}\]
In view of (\ref{sg}) the last form  is equivalent to iii) of
(\ref{twistdef}). 
This proves that $\pi_L\ci S^{-1} \ci S^{\prime}$ is a twist relating
$S^{\prime}$ to $S$.
\qed

Let ${\cal A}=({\cal A}_L,{\cal A}_R,S)$ and ${\cal A}^{\prime}=({\cal
A}_{L^{\prime}}^{\prime},{\cal A}_{R^{\prime}}^{\prime},S^{\prime})$
be Hopf algebroids such 
that the underlying left bialgebroids $ {\cal A}_L$ and ${\cal
A}_{L^{\prime}}^{\prime} $ are isomorphic via the isomorphism $(\Phi:A\to
A^{\prime} , \phi:L\to L^{\prime})$. Then by Proposition \ref{Sisom} also
the underlying right bialgebroids ${\cal A}_R$ and ${\cal
A}_{R^{\prime}}^{\prime}$ are isomorphic and by Theorem \ref{twist}
$S^{\prime}(a^{\prime})= \Phi\ci S\left(\Phi^{-1}(a^{\prime})\lu
g_*\right)$ for a unique twist $g_*$ of $({\cal A}_L,S)$  and all
$a^{\prime}\in {A}^{\prime}$. The Hopf algebroids ${\cal A}$ and
${\cal A}^{\prime}$ are called {\em bialgebroid isomorphic} in the
following. 

\medskip

Recall from \cite{BKSz} that a weak Hopf algebra ${\bf
H}=(H,\Delta,\varepsilon,S)$ over a commutative ring $k$ with
bijective antipode $S$ determines a
Hopf algebroid ${\cal H}=({\cal H}_L,{\cal H}_R,S)$ as follows:
\bea {\cal H}_L&=&(H,L,\id_L,S^{-1}\vert_L,p_L\ci \Delta,\sqcap^L)\lb{hl}\\
     {\cal H}_R&=&(H,R,\id_R,S^{-1}\vert_R,p_R\ci
\Delta,\sqcap^R)\lb{hr}
\eea
where $\sqcap^L:H\to H$ is defined as $h\mapsto \varepsilon(1_{[1]} h)
1_{[2]}$, $1_{[1]}\ot 1_{[2]}=\Delta(1)$, $L\colon = \sqcap^L(H)$, and
$p_L$ is the canonical projection $H\ot_k H\to H\ot_L H$ the $L-L$
-bimodule structure on $H$ being given by 
\[ l\cdot h\cdot l^{\prime}\colon = l S^{-1}(l^{\prime})h.\]
Similarly, $\sqcap^R:H\to H$ is defined as $h\mapsto 1_{[1]} \varepsilon(h
1_{[2]})$, $R\colon = \sqcap^R(H)$, and
$p_R$ is the canonical projection $H\ot_k H\to H\ot_R H$ the $R-R$
-bimodule structure on $H$ being given by 
\[ r\cdot h\cdot r^{\prime}\colon = h r^{\prime} S^{-1}(r).\]
In the view of Theorem \ref{twist} we can obtain examples of Hopf
algebroids by twisting weak Hopf algebras. The twists of the datum
$({\cal H}_L,S)$ are the characters on the weak Hopf algebra ${\bf H}$
or -- if $H$ is finite dimensional as a $k$-space -- the group-like
elements \cite{BNSz,V} in the $k$-dual weak Hopf algebra ${\hat {\bf H}}$.
%This implies that the Connes-Moscovici algebra  ${\cal H}_{FM}$ in
%\cite{ConnMo2} -- which is a twisted weak Hopf algebra -- is a Hopf
%algebroid in the sense of Definition \ref{haid}.

The twistings of cocommutative Hopf algebras are of particular interest
as they provide examples of Hopf algebroids in the sense of Definition
\ref{haid} that do {\em not} satisfy the Lu-Hopf algebroid axioms of
\cite{Lu}. (For \cite{Lu}'s definition see the Introduction above.)
If ${\cal H}_L$ is the left bialgebroid (\ref{hl}) corresponding to the {\em
cocommutative $k$-Hopf
algebra} ${\bf H}$
%$(H,\Delta,\varepsilon,S)$ 
that is ${\cal H}_L=(H,L\equiv
k,s_L\equiv \eta,t_L\equiv \eta,\gamma_L\equiv \Delta,\pi_L\equiv
\varepsilon)$ -- where $\eta:k\to H$ is the unit map $\lambda\mapsto
\lambda 1_H$ -- then the base ring is $k$ itself, so the canonical
projection is the identity map $p_L=\id_{H\ot_k H}$. Then also
$\xi=\id_{H\ot_k H}$. Let $\chi$ be a character on ${\bf H}$ that is an
algebra homomorphism $H\to k$, and $S_{\chi}\colon = (\chi\ot
S)\ci \Delta$ the twisted antipode. Then denoting $\Delta(h)=h_{(1)} \ot
h_{(2)}$ we have
\bea s_L\ci \pi_L(h)&=& \varepsilon(h) 1_H \quad {\rm and}\nn
h_{(1)} S_{\chi} (h_{(2)})&=& 
h_{(1)} \chi(h_{(2)})S(h_{(3)})=
\chi(h_{(1)})h_{(2)} S(h_{(3)})=
\chi(h_{(1)})\varepsilon(h_{(2)})1_H=
\chi(h)1_H
\nonumber\eea
hence $({\cal H}_L,S_{\chi})$ -- which satisfies the conditions
(\ref{lui}-\ref{luiv}) by Theorem \ref{twist}, hence defines a
Hopf algebroid in the 
sense of Definition \ref{haid} -- is a Lu-Hopf algebroid (with the
only possible section $\xi=\id_{H\ot_k H}$) if and only if
$\chi=\varepsilon$. It is easy, however, to find non-trivial
characters on cocommutative Hopf algebras. In the simplest case of the
group Hopf algebra $kZ_2$ we can set $\chi(t)=-t$ (where $t$ is the
second order generator of $Z_2$) if the characteristic of $k$ is
different from $2$. This gives the Example presented in \cite{BKSz}
proving that the two definitions -- the one in \cite{BKSz} and the one
in \cite{Lu} -- of Hopf algebroid are not equivalent.

\medskip

We know from the above considerations that the weak Hopf algebras with
bijective antipodes
provide examples of Hopf algebroids. In the following we identify the
Hopf algebroids that arise this way.

It is proven in \cite{Sch2,Sz} that a left bialgebroid ${\cal
A}_L=(A,L,s_L,t_L,\gamma_L,\pi_L)$ has a weak bialgebra structure {\em
if and only if} $A$ is an algebra over a commutative ring $k$ and $L$
is a  separable $k$-algebra. Indeed let us fix a separability
structure, that is a datum
$(L,k,\delta: L\to L\ot_k L, \psi:L\to k)$ where $\delta$ is a
coassociative coproduct with counit $\psi$ satisfying
\bea (\id_L \ot_k m_L)\ci (\delta\ot_k \id_L)=\delta&\ci& m_L= (m_L\ot_k
\id_L)\ci (\id_L\ot_k \delta)\nn
m_L\ci \delta&=&\id_L.\nonumber \eea
The $m_L$ denotes the multiplication in $L$. The  weak bialgebra
structure on $A$, corresponding to the given separability structure,
reads as
\be \Delta(a)\colon = t_L(e_i)a\di\ot s_L(f_i)a\dii , \qquad
\varepsilon(a)= \psi \ci \pi_L(a)\lb{wba}\ee
where 
$\sum_i e_i\ot f_i=\delta(1_L)$  and the summation symbol is omitted. 

On the other hand if ${\bf H}=(H,\Delta,\varepsilon)$ is a weak 
bialgebra over the commutative ring $k$ then a separability structure
$(L,k,\delta,\varepsilon)$ is given by
\be \delta(l)=l \sqcap^L(1_{[1]}) \ot 1_{[2]}\equiv
                \sqcap^L(1_{[1]}) \ot 1_{[2]} l. \lb{sep}\ee

This implies that the separability of the base algebra $L$ is a
necessary condition for the Hopf algebroid ${\cal A}=({\cal A}_L,{\cal
A}_R,S)$ to have a weak Hopf algebra structure. The different
separability structures determine however different weak bialgebra
structures (\ref{wba}). The given antipode $S$ of the Hopf algebroid ${\cal A}$
cannot make all of them into a weak Hopf algebra. The following
Theorem \ref{wha} gives a criterion on the separability structure
$(L,k,\delta,\psi)$ under which the corresponding weak bialgebra
(\ref{wba}) together with (a twist of) the antipode $S$ becomes a weak
Hopf algebra.

Let  ${\cal A}_L=(A,L,s_L,t_L,\gamma_L,\pi_L)$ be a left bialgebroid such
that  $A$ is an algebra over a commutative ring $k$ and $L$ is a
separable $k$-algebra. Let us fix a separability structure
$(L,k,\delta,\psi)$ and introduce the notation $\delta(1_L)=e_i\ot
f_i$ (summation on $i$ is understood). Let $\Delta$ and $\varepsilon$
be as in (\ref{wba}). Then we can equip the $k$-space ${\hat A}$ of
$k$-linear maps $A\to k$ with an algebra structure with the
multiplication 
\be \phi\phi^{\prime}\colon = (\phi\ot \phi^{\prime})\ci \Delta
\lb{ahat}\ee 
and unit $\varepsilon$. Now we are ready to formulate
\bt  \lb{wha}
Let  ${\cal A}=({\cal A}_L,{\cal A}_R,S)$ be a Hopf algebroid
such that the total ring $A$ is an algebra over a commutative ring $k$
and the base ring $L$ of ${\cal A}_L$ is a separable $k$-algebra. Then
fixing a separability structure $(L,k,\delta,\psi)$ the corresponding
weak bialgebra (\ref{wba}) and the antipode $S$ form a weak Hopf
algebra if and only if $\psi\ci\pi_L\ci S=\psi\ci\pi_L$. Furthermore,
there exists a twisted antipode $S_g$ making the weak bialgebra
(\ref{wba}) a weak Hopf algebra if and only if the element
$\psi\ci\pi_L\ci S$ of the algebra ${\hat A}$ -- defined in
(\ref{ahat}) -- is invertible.
\et

\pr {\em if part:} The separability structure $(L,k,\delta,\psi)$
defines an isomorphism of the algebras $\atld$ 
-- the $L$-dual algebra (\ref{ldual}) of ${\cal A}_L$ --
and ${\hat A}$ in (\ref{ahat}): 
\[\begin{array}{rrcl}
\kappa:& {\hat A}\ \to& \atld \qquad 
\phi &\mapsto\  \left( a\ \mapsto\ \phi(t_L(e_i)a)f_i \right)\nn
\kappa^{-1}:& \atld \ \to& {\hat A} \qquad  \phi_*\!&\mapsto\  \psi\ci
\phi_*.\end{array}\]
The element $\kappa(\psi\ci \pi_L\ci S)$ satisfies the properties {\em
i)-iii)} of Definition \ref{twistdef}:
\bea 1_A&\lu& \kappa(\psi\ci\pi_L\ci S)= s_L(f_i)\psi\ci \pi_L\ci S\ci
t_L(e_i)=s_L(1_L)=1_A\nn 
(a&\lu&\kappa(\psi\ci\pi_L\ci S))(b\ \lu\ \kappa(\psi\ci\pi_L\ci S))= \nn
&=&s_L(f_i)a\dii s_L(f_j) b\dii \psi\ci \pi_L\ci S(t_L(e_i)a\di)\psi\ci
\pi_L\ci S(t_L(e_j)b\di) =\nn
&=&s_L(f_i)a\dii b\dii \psi\ci \pi_L \ci S \left(t_L[e_j
\psi(f_j\pi_L\ci S\{t_L(e_i)a\di\})]b\di\right)=\nn
&=&s_L(f_i)a\dii b\dii \psi\ci \pi_L\ci S(t_L(e_i)a\di b\di)=ab\ \lu\ 
\kappa(\psi\ci\pi_L\ci S)\nn 
S(a\di)&\lu& \kappa(\psi\ci\pi_L\ci S) \ot a\dii \ \lu\ 
\kappa(\psi\ci\pi_L\ci S) =\nn 
&=&s_L(f_i)S(a\di)\dii \ot s_L(f_j){a\dii}\dii 
\psi\ci \pi_L\ci S\left(t_L(e_i)S(a\di)\di\right)
\psi\ci \pi_L\ci S\left(t_L(e_j){a\dii}\di\right)=\nn
%&=&s_L(f_i)S({a\di} \ui)\ot {a\dii}\dii 
%\psi\ci \pi_L\ci S\left(t_L[e_j \psi(f_j \pi_L\ci S(t_L(e_i)S({a\di}
%\uii)))]{a\dii}\di \right)=\nn
%&=& s_L(f_i)S({a\di} \ui)\ot {a\dii}\dii 
%\psi\ci \pi_L\ci S\left(s_R\ci \pi_R\ci S[{a\di}\uii S^{-2}\ci
%t_L(e_i)]{a\dii}\di\right)= \nn
%&=& s_L(f_i)S({a} \ui)\ot {{a\uii}\dii}\dii 
%\psi\ci \pi_L\ci S\left(s_R\ci \pi_R\ci S[{a\uii}\di S^{-2}\ci
%t_L(e_i)]{{a\uii}\dii}\di\right)= \nn
%&=& s_L(f_i)S({a} \ui)\ot {{a\uii}\dii}\dii 
%\psi\ci \pi_L\left(t_L(e_i)S( {{a\uii}\dii}\dii)s_L\ci \pi_L\ci
%S^2({a\uii})\di\right)=\nn 
%&=& \psi\left( \pi_L\ci S[s_R\ci \pi_R\ci
%S({a\uii}\di){{a\uii}\dii}\di]e_i\right) s_L(f_i) S(a\ui)\ot
%{{a\uii}\dii}\dii=\nn
%&=&S\left(a\ui s_R\ci\pi_R[S({{a\uii}\di}\dii){{a\uii}\di}\dii]\right)\ot
%{a\uii}\dii= 
%%%%%%%%%%%%%%%%%%%%%%%%%%%%%%%%%%%%%%%%%%%%%%%%%%%%%%%%%%%%%%%%%%%%%%%%%
&=&s_L(f_i)S(a\di)\dii  s_L(f_j)\ot {a\dii}\dii 
\psi\ci \pi_L\ci S\left(t_L(e_i)S(a\di)\di\right)
\psi\ci \pi_L\ci S\left(t_L(e_j){a\dii}\di\right)=\nn
&=&s_L(f_i)S({a\di}\ui) \ot {a\dii}\dii 
\psi\left( f_j \pi_L\ci S\left[t_L(e_i)S({a\di}\uii)\right]\right)
\psi\ci \pi_L\ci S\left(t_L(e_j){a\dii}\di\right)=\nn
&=&s_L(f_i)S({a\di}\ui) \ot {a\dii}\dii 
\psi\ci \pi_L\ci S\left(t_L(e_i)S({a\di}\uii){a\dii}\di\right)=\nn
&=&s_L(f_i)S({a}\ui) \ot {a\uii}\dii 
\psi\ci \pi_L\ci
S\left(t_L(e_i)S({{a\uii}\di}\di){{a\uii}\di}\dii\right)=\nn 
&=&s_L(f_i)S({a}\ui) \ot {a\uii}\dii 
\psi\ci\pi_L\ci S\ci t_L\left(\pi_L\ci S({{a\uii}\di})e_i\right)=\nn 
&=&\psi\left(\pi_L\ci S({{a\di}\uii})e_i\right)s_L(f_i)S({a\di}\ui) 
\ot a\dii=\nn
&=&S\ci s_R\ci \pi_R({a\di}\uii)S({a\di}\ui)\ot a\dii=
S(a\di)\ot a\dii
\nonumber \eea
hence -- since $\psi\ci \pi_L\ci
S\in {\hat A}$ is invertible by assumption -- both  $\kappa(\psi\ci
\pi_L\ci S)$ and its inverse $g_*\colon =  \kappa(\psi\ci \pi_L\ci
S)^{-1}$ are twists of $({\cal A}_L,S)$ in the sense of Definition
\ref{twistdef}. 

Using the standard properties \cite{W} of the quasi-basis $e_i\ot f_i$
of $\psi$ one checks that the twisted antipode $S_g:a\mapsto S(a\lu
g_*)$ makes the weak bialgebra (\ref{wba}) a weak Hopf algebra. That
is by definition 
\[ \psi\ci \pi_L\ci S_g=
\kappa^{-1}(g_*)(\psi\ci\pi_L\ci S)=\psi\ci \pi_L.\]
Since $S_g\ci t_L=S\ci t_L=s_L$ we have
\[ \psi\left(\pi_L(a) \pi_L\ci S_g^2\ci t_L(l)\right)=
\psi\ci\pi_L\left( S_g^2\ci t_L(l) s_L\ci \pi_L(a)\right)=
\psi\ci\pi_L\left(t_L\ci \pi_L(a) s_L(l)\right)=
\psi(l\pi_L(a)).\]
This implies that $e_i\ot f_i=f_i\ot \pi_L\ci S_g^{-2}\ci t_L(e_i)$ and
hence $e_i\ot s_L(f_i)=f_i\ot S_g^{-1}\ci t_L(e_i)$. Then
\bea m_A\ci (S_g\ot_k \id_A)&\ci& \Delta(a) =
S\left( s_L\ci g_*(a\di)t_L(e_i)a\dii\right)s_L(f_i)a_{(3)}=
s_R\ci \pi_R(a\lu g_*)=\nn
&=&t_L\ci\pi_L\ci S_g(a)=
\psi\left(\pi_L\ci S_g(a)e_i\right)t_L(f_i)=\nn
&=&\psi\ci \pi_L\left( a S_g^{-1}\ci t_L(e_i)\right) t_L(f_i)=
t_L(e_i)\psi\ci \pi_L(as_L(f_i))\nn
m_A\ci (\id_A \ot_k S_g)&\ci& \Delta(a)=
t_L(e_i)a\di S(a_{(3)})S\ci s_L\ci g_*(s_L(f_i)a\dii)=\nn
&=&t_L(e_i)a\di t_L(e_k)\psi\ci \pi_L\ci S\left({a_{(3)}}\ui
t_L(f_k)\right)
S\left(s_L\ci g_*(s_L(f_i)a\dii) {a_{(3)}}\uii\right)=\nn
&=&t_L(e_i)a\di t_L(e_k)\psi\ci \pi_L\ci S\left( 
s_L\ci g_*[s_L(f_i){{a\dii}\ui}\di] {{a\dii}\ui}\dii
t_L(f_k)\right)S({a\dii}\uii)=\nn
&=&t_L(e_i)a\di t_L(e_k)\psi\ci \pi_L\left(s_L(f_i) {a\dii}\ui
t_L(f_k)\right) S({a\dii}\uii)=\nn
&=&t_L\ci \pi_L({a\ui}\dii t_L(f_k)){a\ui}\di t_L(e_k)S(a\uii)=
s_L\ci \pi_L(a)=
\psi\ci\pi_L(t_L(e_i)a)s_L(f_i)\nn
S_g(t_L(e_i)a\di)&s_L&(f_i) t_L(e_j)a\dii S_g(s_L(f_j)a_{(3)})=
S_g(a\di)s_L\ci \pi_L(a\dii)=S_g(a)
\nonumber\eea
for all $a\in A$. This finishes the proof of the if part.

{\em only if part:}
If ${\cal H}_L$ is a left bialgebroid (\ref{hl}) corresponding to the
weak Hopf algebra ${\bf H}$  then its
base ring $L$
is a separable $k$-algebra \cite{BNSz,Sch2}. The separability structure
%$(L,k,\delta,\psi)$ can be given in terms of the weak Hopf data:
%$\delta(l)\colon = lS(1_{[1]})\ot_k 1_{[2]}\equiv S(1_{[1]})\ot_k
%1_{[2]}l$ and $\psi(l)=\varepsilon (l)$ for all $l\in L$ -- where
%$1_{[1]}\ot_k 1_{[2]}=\Delta(1)$. 
(\ref{sep}) is determined by the weak Hopf data in ${\bf H}$. Now if the datum 
$({\cal H}_L,S^{\prime})$ is obtained as a twist of the datum $({\cal
H}_L,S)$ then the twist element relating them is of the form
\[\sqcap^L\ci S^{-1}\ci S^{\prime}=
\kappa(\varepsilon\ci\sqcap^L\ci S^{-1}\ci S^{\prime})= 
\kappa(\varepsilon\ci S^{-1}\ci S^{\prime})= 
\kappa(\varepsilon\ci   S^{\prime})= 
\kappa(\varepsilon\ci\sqcap^L\ci S^{\prime})\]
which is an invertible element of ${\cal H}_*$, by definition. Since
$\kappa$ is an algebra isomorphism this proves that
$\varepsilon\ci\sqcap^L\ci 
S^{\prime}$ is an invertible element of ${\hat H}$, hence the claim.
\qed

\br Proposition \ref{wha} implies that the Hopf algebroid ${\cal
A}=({\cal A}_L,{\cal A}_R,S)$ is obtained by a twist of a {\em Hopf
algebra} over the commutative ring $k$ if and only if the following
conditions hold true:
\bea &i)&L\simeq R^{op} \ {\rm is\ isomorphic\ to \ } k\nn
     &ii)&A\ {\rm is\ a\ } k-{\rm algebra}\nn
     &iii)&\pi_L\ci S:A\to k \ {\rm is\ invertible\ in\ } {\hat A}, {\rm \
the \ } k-{\rm dual\ algebra\ of \ } A
\nonumber \eea
In particular it is a Hopf algebra if and only if the conditions $i)$
and $ii)$ hold true and $\pi_L\ci S=\pi_L$.
\er
 
\section{The theory of integrals}
\lb{int}
\setc{0}

The left/right {\em integrals} in a Hopf algebra $(H,\Delta,\varepsilon,S)$
are the invariants of left/right regular module. That is $\ell/\err\in H$ is
a left/right  integral if 
\[ h\ell=\varepsilon (h) \ell\quad /\quad \err h =\err\varepsilon (h)\]
for all $h\in H$. This notion has been generalized to a weak Hopf
algebra $(H,\Delta,\varepsilon ,S)$ as $\ell/\err\in H$ is a
left/right  integral if 
\[ h\ell=\sqcap^L(h) \ell\quad /\quad \err h=\err\sqcap^R(h)\]
for all $h\in H$ where $\sqcap^L(h)=\varepsilon(1_{[1]} h)1_{[2]}$,
$\sqcap^R(h)=1_{[1]} \varepsilon(h1_{[2]})$ and $1_{[1]}\ot
1_{[2]}=\Delta(1_H)$.  

It is then straightforward to generalize the notion of integrals to
Hopf algebroids:
\bd The left integrals in a left bialgebroid ${\cal
A}_L=(A,L,s_L,t_L,\gamma_L,\pi_L)$ are the elements $\ell\in A$ that
satisfy
\[ a\ell=s_L\ci \pi_L(a) \ell \]
for all $a\in A$. The right ideal  of left integrals is denoted by
$\ila$. The right integrals in the right bialgebroid ${\cal 
A}_R=(A,R,s_R,t_R,\gamma_R,\pi_R)$ are the elements $\err\in A$ for which
\[ \err a=\err s_R\ci \pi_R(a )\]
for all $a\in A$. The left ideal of right integrals is denoted by ${\cal
I}^R(A)$. 

In a Hopf algebroid ${\cal A}=({\cal A}_L,{\cal A}_R,S)$ the left/right
integrals are the left/right integrals in ${\cal A}_L$/${\cal A}_R$.
\ed

As a support of this definition the Lemma 3.2 of \cite{BNSz} generalizes as
\bl \lb{intpr}
The following are equivalent:
\[ \begin{array}{rl}
i)&\ell\in  {\cal I}^L({\cal A})\nn
ii)&a\ell = t_L\circ \pi_L(a)\ell \quad {for\ all\ } a\in A\nn
iii)& S(\ell)\in {\cal I}^R({\cal A})\nn
iv)& S^{-1} (\ell)\in {\cal I}^R({\cal A})\nn
v) & S(a) \ell^{(1)} \otimes \ell^{(2)} = \ell^{(1)} \otimes a \ell^{(2)} 
{\ as\  elements\  of\ } A^R \otimes ^RA {\ for\ all\ } a\in A.
\nonumber\end{array}\]
\el
A left integral in a Hopf algebroid ${\cal A}$ is also a left integral
in ${\cal A}_{cop}$ and it is a right integral in ${\cal A}^{op}$.

For the Hopf algebroid ${\cal A}=({\cal A}_L,{\cal A}_R,S)$ we
introduce the following notation: Let $\asrd$ and $\atrd$ denote the
dual rings (\ref{rdual}) of the right bialgebroid ${\cal A}_R$ and
$\asld$ and $\atld$ denote the dual rings (\ref{ldual}) of the left
bialgebroid ${\cal A}_L$. 
We define the non-degeneracy of an integral as follows:
\bd \lb{nddef}
The left integral $\ell \in \ila$ is {\em non-degenerate} if the
maps 
\bea \fsr: \asrd\ \to &A\qquad \phisr&\mapsto\ \phisr\ru \ell\quad {\rm
and}\lb{fsr}\\
\ftr:\  \atrd\ \to &A \qquad \phitr &\mapsto \ \phitr\rd \ell
\lb{ftr}\eea 
are bijective. A right integral $\err\in {\cal I}^R({\cal A})$ is
non-degenerate if the maps
\bea {_L\err}: \asld\to &A\qquad \phisl&\mapsto\ \err\ld \phisl\quad
{\rm and}\\
     {\err_L}: \atld\to &A\qquad \phitl&\mapsto\ \err\lu \phitl
\eea
are bijective.
\ed
The identities 
\[ \begin{array}{rcl}
S(\phisr\ru a)\ =&S(a)\ld\pi_L\ci s_R\ci \phisr\ci S^{-1}\qquad
S^{-1}(\phisr\ru a)&=\ S^{-1}(a)\ld \pi_L\ci t_R\ci \phisr\ci S\\
     S(\phitr\rd a)\ =&S(a)\lu \pi_L\ci s_R\ci \phitr\ci S^{-1}\qquad
S^{-1}(\phitr\rd a)&=\ S^{-1}(a)\lu \pi_L\ci t_R\ci \phitr\ci S
\end{array}\]
imply that $\ell$ is a non-degenerate left integral if and only if
$S(\ell)$ is a non-degenerate right integral and if and only if
$S^{-1}(\ell)$ is a non-degenerate right integral.

Let $\ell$ be a non-degenerate left integral in the Hopf algebroid
${\cal A}$. Introducing $\lsr\colon = \fsr^{-1}(1_A)$ and $\ltr\colon
=\ftr^{-1}(1_A)$ we have
\bea (\lsr\lu S(a))\ru \ell&=& 
\ell\uii t_R\ci \lsr(S(a)\ell\ui)=
a(\lsr \ru \ell)=a\nn
(\ltr\ld S^{-1}(a))\rd \ell&=&
\ell\ui s_R\ci \ltr(S^{-1}(a)\ell\uii)=
a(\ltr\rd \ell)=a\nonumber\eea
hence 
\be \fsr^{-1}(a)=\lsr\lu S(a)\qquad 
    \ftr^{-1}(a)=\ltr\ld S^{-1}(a).\lb{fsrinv}\ee

Recall that in the Definition \ref{nddef} the non-degeneracy of a {\em
left} integral is defined in terms of the duals $\asrd$ and $\atrd$ of
the {\em right} 
bialgebroid ${\cal A}_R$. The explanation of this is that -- in the view
of property $v)$ in Lemma \ref{intpr} -- {\em this} notion
of non-degeneracy is equivalent to the property that $(\lsr,\ell\ui\ot
S(\ell\uii))$ is a Frobenius system for the ring extension $s_R:R\to
A$. This  implies in particular that for a Hopf algebroid ${\cal A}$
possessing a non-degenerate left integral $\ell$ the
modules $\asr$, $\atr$, $\asl$ and $\atl$ are all finitely
generated projective. Hence by the results in \cite{KSz} the
corresponding duals $\asrd$ and $\atrd$ carry left- and $\asld$ and
$\atld$ carry right bialgebroid structures. In addition to the maps
$\fsr$ and $\ftr$ also the maps
\[\begin{array}{cccc}
 \ftl: \atld \to A  &\qquad \phitl&\mapsto& \quad\ell\leftharpoonup
\phitl \nn
\fsl: \asld \to A  &\qquad   \phisl&\mapsto& \quad \ell\leftharpoondown
\phisl
\end{array}\]
turn out to be bijective. The map 
\be {\tilde S}:a\ \mapsto\ 
%S^{-1}\left((\ltr \ld \ell)\rd a\right)
\ell\lu(a\ru \pi_L\ci s_R\ci \lsr) \lb{twap}\ee
is an anti-automorphism of the
ring $A$ and we have the following isomorphisms of left bialgebroids:

\smallskip

\begin{center}
\begin{picture}(252,80)
\put(60,0){$\asrd_L$}
\put(214,0){$\atrd_L$}
\put(59,68){$(\atld_R)^{op}_{cop}$}
\put(207,68){$(\asld_R)^{op}_{cop}$}
\put(79,2){\vector(1,0){130}}
\put(92,70){\vector(1,0){115}}
\put(70,65){\vector(0,-1){55}}
\put(217,65){\vector(0,-1){55}}
\put(-20,32){$(\fsr^{-1}\ci\ftl,\pi_R\ci s_L)$}
\put(222,32){$(\ftr^{-1}\ci \fsl,\pi_R\ci t_L)$}
\put(73,8){$(\ftr^{-1}\ci {\tilde S}^{-1}\ci \fsr,\pi_R\ci S^{-1}\ci t_R)$}
\put(100,76){$(\fsl^{-1}\ci {\tilde S}^{-1} \ci \ftl,\id_R)$}
\end{picture}
\end{center}

\bigskip

Let ${\cal A}=({\cal A}_L,{\cal A}_R,S)$ be a Hopf algebroid with a
non-degenerate left integral $\ell$ and let ${\cal A}^{\prime}= ({\cal
A}_{L^{\prime}}^{\prime},{\cal A}_{R^{\prime}}^{\prime},S^{\prime})$
be a Hopf algebroid which 
is bialgebroid isomorphic to ${\cal A}$ via the isomorphism
$(\Phi:A\to A^{\prime},\phi:L\to L^{\prime})$ of left
bialgebroids. Then $\Phi(\ell)$ is a non-degenerate left integral in
${\cal A}^{\prime}$. 

The antipode $S_*\colon = \ell_L^{-1}\ci {\tilde S}\ci \ell_L$, $\phitl
\ \mapsto\ (\ell\lu \phitl)\ru \pi_L\ci s_R\ci \lsr$ makes the dual
right bialgebroid $\atld_R$ into a Hopf algebroid called ${\cal A}^{\ell}_*$
possessing a two-sided non-degenerate integral $\pi_L\ci s_R\ci \lsr$.
Clearly the bialgebroid isomorphism class of ${\cal A}^{\ell}_*$ does
not depend on the choice of the non-degenerate integral $\ell$. 
The {\em dual of the bialgebroid isomorphism class of ${\cal A}$} is then
defined to be the bialgebroid isomorphism class of ${\cal
A}^{\ell}_*$. This notion of duality is shown to be involutive and
reproduces the duality of finite weak Hopf algebras \cite{BNSz} as
follows: 

Let ${\bf H}=(H,\Delta,\varepsilon,S)$ be a finite weak Hopf algebra over
the commutative ring $k$ and let ${\hat {\bf H}}=({\hat H}, {\hat
\Delta}, {\hat \varepsilon},{\hat S})$ be its $k$-dual weak Hopf algebra
\cite{BNSz}. The corresponding Hopf algebroids are denoted by ${\cal
H}=({\cal H}_L,{\cal H}_R,S)$ and ${\hat {\cal H}}=({\hat {\cal
H}}_{\hat L},{\hat {\cal H}}_{\hat R},{\hat S})$, respectively. Then
the right bialgebroids 
${\cal H}_{*R}$ and ${\hat {\cal H}}_{\hat R}$ are isomorphic via
\bea \Phi:{\cal H}_*\to& {\hat H}\qquad 
\psi_*&\mapsto \varepsilon \ci \psi_*\lb{whadiso}\\
\phi: L\ \ \to & \!\!\!{\hat R} \ \qquad
l &\mapsto \varepsilon_{[1]}\ \varepsilon_{[2]}(l) \eea
where $\varepsilon_{[1]}\ot \varepsilon_{[2]}={\hat
\Delta}(\varepsilon)$. This implies that the Hopf algebroids ${\cal
H}_*^{\ell}$ and ${\hat {\cal H}}$ are bialgebroid isomorphic for any
choice of the non-degenerate left integral $\ell$. Making use of the
separability structure (\ref{sep}) we can equip 
${\cal H}_*$ with a weak bialgebra structure. Since the weak Hopf
structure on a weak bialgebra is unique, there is a unique twist of $S_*$
that makes this weak bialgebra into a weak Hopf algebra. The map  $\Phi$
in (\ref{whadiso}) is a weak Hopf algebra isomorphism from it to
${\hat{\bf H}}$.

\section{Antipode from non-degenerate integral}
\lb{LS}
\setc{0}

The Larson-Sweedler theorem \cite{LaSw} on Hopf algebras
states that a finite dimensional bialgebra over a field 
is a Hopf algebra if and
only if it has a non-degenerate left integral. It has been generalized
to weak Hopf algebras in \cite{V}. At the moment no generalization to
Hopf algebroids is known. In this section we present a proof of the
`easy half' of the Larson-Sweedler theorem. Namely we show that if a
finite {\em right} bialgebroid ${\cal A}_R$ -- that is a right
bialgebroid such that $\asr$ and $\atr$ are  finitely generated
projective -- has a 
non-degenerate {\em left} integral then it can be made a Hopf
algebroid. 

Notice that though we called it the `easy half' of the Larson-Sweedler
theorem, it is non-trivial in the sense that left integrals have not
been defined in right bialgebroids until now. Recall that if ${\cal
A}$ is a Hopf algebroid possessing a non-degenerate left integral
$\ell$ then both maps $\fsr$ and $\ftr$ (\ref{fsr}-\ref{ftr}) are
bijective and for all $a\in A$
\[ \ell\ui \ot a\ell \uii= S(a)\ell\ui \ot \ell\uii .\]
The equations (\ref{fsrinv}) imply that the antipode and its inverse
can be written into the forms 
\bea S(a)&=& (\ltr \ld a)\rd \ell\nn
     S^{-1}(a)&=& (\lsr\lu a)\ru \ell
\nonumber\eea
with the help of $\lsr=\fsr^{-1}(1_A)$ and $\ltr=\ftr^{-1}(1_A)$.
These formulae motivate
\bd \lb{bgdnd} An element $\ell$ of a finite right bialgebroid ${\cal
A}_R$ is a 
{\em non-degenerate left integral} if
\bea
&i)& {both\  maps\ } \fsr {\ and\ } \ftr {\ are\  bijective}\\
&ii)& \ell\ui \ot a \ell\uii = 
%\ell^{(1)^{\prime}}\ell\ui \ot t_R\ci \ltr(a \ell^{(2)^{\prime}})
%\ell\uii
[(\ltr\ld a)\rd \ell]\ell\ui\ot \ell\uii\lb{sf}\\
&&a\ell\ui \ot \ell\uii= 
%s_R\ci \lsr (a\ell^{(1)^{\prime}})\ell\ui \ot \ell^{(2)^{\prime}} 
%\ell\uii
\ell\ui \ot [(\lsr\lu a)\ru \ell]\ell\uii \lb{sb}\eea
as elements of $\asr\ot \atr$ for all $a\in A$, where
$\lsr=\fsr^{-1}(1_A)$ and $\ltr=\ftr^{-1}(1_A)$. 
\ed
The following Lemma is of technical use:
\bl \lb{lac} Let ${\cal A}_R=(A,R,s_R,t_R,\gamma_R,\pi_R)$ be a finite
right bialgebroid and $k\in A$ such that the map
\[ k_R:\asrd\to A\qquad \phisr\mapsto\ \phisr\ru k \]
is bijective. Set $\kappa^*\colon = k_R^{-1}(1_A)$. Then
\be  \kappa^*\ru a= s_R\ci \kappa^* (a) \ee
for all $a\in A$. Analogously, for $k\in A$ such that the map
\[ {_R k}:\atrd\to A\qquad \phitr\mapsto\ \phitr\rd k \]
is bijective set ${^* \kappa}\colon = {_R k}^{-1}(1_A)$. Then
\be  {^* \kappa}\rd a= t_R\ci {^* \kappa} (a) \lb{tac}\ee
for all $a\in A$.
\el

\pr Introduce the map ${\hat s}:R\to \asrd$ via ${\hat
s}(r)(a)=r\pi_R(a)$. With its help
\bea \phisr\kappa^*&=& k_R^{-1}( \phisr\kappa^*\ru k)=
k_R^{-1}( \phisr\ru 1_A)=k_R^{-1}( {\hat s}\ci\phisr(1_A) \ru
1_A)=
\nn &=&  
k_R^{-1}({\hat s}\ci\phisr(1_A)  \kappa^*\ru k )={\hat
s}\ci\phisr(1_A)\kappa^* \nonumber
\eea 
for all $\phisr\in \asrd$.  This implies that 
\[ \phisr(\kappa^* \ru a) = (\phisr\kappa^*)(a)=
({\hat s}\ci\phisr(1_A)\kappa^*)(a)=
\phisr(1_A)\kappa^*(a)=\phisr(s_R\ci \kappa^*(a))\]
for all $a\in A$ and $\phisr\in \asrd$. Since $\asr$ is finitely
generated projective by assumption, this proves that $\kappa^*\ru a=s_R\ci
\kappa^*(a)$.  The identity (\ref{tac}) follows by applying the same
proof in ${\cal A}_{cop}$.
\qed

We are ready to formulate
\bt \lb{LS/2}
Let ${\cal A}_R$ be a finite right bialgebroid with non-degenerate
left integral $\ell$. 
Let $\lsr\colon = \fsr^{-1}(1_A)$ and $\ltr\colon = \ftr^{-1}(1_A)$. 
Then the map $S(a)\colon = (\ltr\ld a)\rd \ell$ 
is an antipode making ${\cal A}_R$  into a Hopf algebroid. The
element $\ell$ is a non-degenerate left 
integral in the resulting Hopf algebroid in the sense of Definition
\ref{nddef}. 
\et

\pr We check that $({\cal A}_R,S)$ satisfies the right analogues of
the conditions (\ref{lui}-\ref{luiv}) -- what implies the claim.
The conditions (\ref{sf}) and (\ref{sb}) imply that the map
$S\colon = (\ltr\ld a)\rd \ell$ is bijective with inverse  
$S^{-1}(a)\colon = (\lsr\lu a)\ru \ell$. It is anti-multiplicative by
\[ S(b)S(a)=S(b)S(a)\ell\ui s_R\ci \ltr(\ell\uii)=
\ell\ui s_R\ci \ltr (ab \ell\uii)=
S(ab)\ell\ui s_R\ci \ltr(\ell\uii)=S(ab).\]
Also
\[ S\ci t_R(r)= \ell\ui s_R\ci \ltr(t_R(r)\ell\uii)=s_R(r).\]
Using Lemma \ref{lac}
\[ a\ui S(a\uii)=
a\ui \ell\ui s_R\ci \ltr(a\uii\ell\uii)=
\ltr \rd a\ell=
t_R\ci \ltr (a\ell)=
t_R\ci \pi_R\left((\ltr\ld a)\rd \ell \right)=
t_R\ci \pi_R\ci S(a).\]
%This implies that 
%\[ a\ell =a\ui \ell\ui s_R\ci \pi_R (a\uii \ell\uii)=a\ui
%S(a\uii)\ell=
%t_R\ci \pi_R\ci S(a) \ell\]
Using the identities
\bea \gamma_R\ci S(a)&=&\ell\ui\ot (\ltr\ld a)\rd \ell\uii\lb{grs}\\
\gamma_R\ci S^{-1}(a)&=& (\lsr \lu a)\ru \ell\ui \ot \ell\uii\lb{grsi}
\eea
one shows that
\bea S\left(S^{-1}(a)\uii \right)&\ot&  S\left(S^{-1}(a)\ui \right)\ui 
\ot  S\left(S^{-1}(a)\ui \right)\uii=\nn
%&=&S(\ell\uii)\ot S\left((\lsr\lu a)\ru \ell\ui\right)\ui
%\ot S(\ell\uii)\ot S\left((\lsr\lu a)\ru \ell\ui\right)\uii=\nn
%&=&
%S(\ell\uii)\ot \ell^{(1)^{\prime}}\ot \left(\ltr\ld[(\lsr\lu a)\ru
%\ell\ui]\right)\rd \ell^{(2)^{\prime}}=\nn
&=&
S(\ell^{(3)})\ot  \ell^{(1^{\prime})}\ot \ell^{(2^{\prime})}
s_R\ci \ltr \left( \ell\uii t_R\ci \lsr
(a\ell\ui)\ell^{(3^{\prime})}\right)=\nn
%&=&S(\ell^{(3)})\ot  \ell^{(1)^{\prime}}\ot 
%s_R\ci \lsr (a\ell\ui) \ell^{(2)^{\prime}}s_R\ci \ltr \left(
%\ell\uii\ell^{(3)^{\prime}}\right)=\nn 
&=&S(\ell^{(3)})\ot  \ell^{(1^{\prime})}\ot (\lsr\ru a\ell\ui)
 \ell^{(2^{\prime})}s_R\ci \ltr \left(
\ell\uii\ell^{(3^{\prime})}\right)=\nn  
%&=&S(\ell^{(4)})\ot  \ell^{(1)^{\prime}}\ot
%a\uii \ell\uii t_R\ci \lsr(a\ui\ell\ui)  \ell^{(2)^{\prime}}s_R\ci
%\ltr \left( \ell^{(3)} \ell^{(3)^{\prime}}\right)=\nn   
&=&S(\ell^{(3)})\ot  \ell^{(1^{\prime})}\ot
a\uii \left( \ltr\rd \ell\uii t_R\ci \lsr(a\ui\ell\ui)
\ell^{(2^{\prime})} \right)=\nn
%&=&S(\ell^{(3)})\ot  \ell^{(1)^{\prime}}\ot
%a\uii t_R\ci \ltr \left( \ell\uii t_R\ci \lsr(a\ui\ell\ui)
%\ell^{(2)^{\prime}}\right)=\nn
&=&S(\ell^{(3)})\ot   \ell^{(1^{\prime})} s_R\ci \ltr \left( \ell\uii
t_R\ci \lsr(a\ui\ell\ui) \ell^{(2^{\prime})}\right) \ot a\uii=\nn
%&=&S(\ell\uii)\ot \left( \ltr \ld [ ( \lsr \lu a\uii)\ru
%\ell\ui]\right)\rd \ell  \ot a\uii=\nn
%&=&S(\ell\uii)\ot S\left(  ( \lsr \lu a\uii)\ru \ell\ui \right )  \ot
%a\uii=\nn
&=&S\left(S^{-1}(a\ui)\uii\right) \ot S\left(S^{-1}(a\ui)\ui\right)
 \ot a\uii
\lb{coa1}\eea
where $\ell^{(1^{\prime})}\ot
\ell^{(2^{\prime})}=\gamma_R(\ell)=\ell\ui\ot \ell\uii$. Hence
\be S^2\left(S^{-1}(a)\uii\right)S\left(S^{-1}(a)\ui\right)\ui \ot
S\left(S^{-1}(a)\ui\right)\uii = s_R\ci \pi_R(a\ui)\ot a\uii.\lb{lem}\ee
Using (\ref{grs}), (\ref{grsi}) and (\ref{lem}) one checks that
\bea 
S\left( S^{-1}(a)\uii \right) &\ot&   S\left( S^{-1}(a)\ui \right) =\nn
&=&S\left( S^{-1}(a)\uii \right) \ot 
S^{-1} \left( S^2( S^{-1}(a)\ui)\ui s_R\ci \pi_R 
[ S^2(S^{-1}(a)\ui)\uii]\right)=\nn
%&=&S\left( S^{-1}(a)\uii \right) \ot 
%t_R\ci \pi_R [ S^2(S^{-1}(a)\ui)\uii] S^{-1} 
%\left( S^2( S^{-1}(a)\ui)\ui\right) =\nn
&=&S^{-1}\ci t_R\ci \pi_R [ S^2(S^{-1}(a)\ui)\uii]
S\left( S^{-1}(a)\uii \right) \ot   
S^{-1} \left( S^2( S^{-1}(a)\ui)\ui\right) =\nn
%&=&S^{-1}\ci t_R\ci \pi_R[( \ltr\ld  S(S^{-1}(a)\ui))\rd \ell\uii]
%S \left( S^{-1}(a)\uii \right) \ot  S^{-1}(\ell\ui)=\nn
&=&S^{-1}\ci t_R\ci \ltr \left( S(S^{-1}(a)\ui)\ell\uii\right)
S \left(S^{-1}(a)\uii \right)  \ot  S^{-1}(\ell\ui)=\nn
%&=&S^{-1}\left( \ltr\rd  S(S^{-1}(a)\ui)\ell\uii\right)
%S \left(S^{-1}(a)\uii \right)  \ot  S^{-1}(\ell\ui)=\nn
%S^{-1}\left(S(S^{-1}(a)\ui)\ui \ell\uii 
%s_R\ci \ltr[S(S^{-1}(a)\ui)\uii\ell^{(3)} ] \right)
%S \left(S^{-1}(a)\uii \right)  \ot  S^{-1}(\ell\ui)=
&=&S^{-1}\left( S^2(S^{-1}(a)\uii)S(S^{-1}(a)\ui)\ui \ell\uii 
s_R\ci \ltr[S(S^{-1}(a)\ui)\uii\ell^{(3)} ] \right)
 \ot  S^{-1}(\ell\ui)=\nn
&=&S^{-1}\left( s_R\ci \pi_R(a\ui)\ell\uii s_R\ci \ltr(a\uii \ell^{(3)})\right)
 \ot  S^{-1}(\ell\ui)=\nn
&=&S^{-1}\left(\ell\uii s_R\ci \ltr(a\ell^{(3)}) \right)  \ot
S^{-1}(\ell\ui)=
S^{-1}\left(S(a)\uii\right)\ot S^{-1}\left(S(a)\ui\right).
\eea
Introducing $\gamma_L(a)\colon = S\left( S^{-1}(a)\uii \right) \ot
S\left( S^{-1}(a)\ui \right) \equiv S^{-1}\left(S(a)\uii\right)\ot
S^{-1}\left(S(a)\ui\right)$ the identity
\[ (\id_A\ot \gamma_R)\ci \gamma_L=(\gamma_L\ot \id_A)\ci \gamma_R\]
is proven by (\ref{coa1}). The last axiom $(\gamma_R\ot \id_A)\ci
\gamma_L= (\id_A\ot \gamma_L)\ci \gamma_R$ is checked similarly:
\bea S^{-1}\left(S(a)\uii\right)\ui &\ot&  S^{-1}\left(S(a)\uii\right)\uii
\ot  S^{-1}\left(S(a)\ui\right)=\nn
%&=& (\lsr\lu S(a)\uii)\ru \ell\ui \ot \ell\uii \ot
%S^{-1}\left(S(a)\ui\right)=\nn 
%&=& \left( \lsr\lu[(\ltr\ld a)\rd \ell^{(2)^{\prime}}]\right)\ru
%\ell\ui \ot  \ell\uii \ot S^{-1}\left(\ell^{(1)^{\prime}}\right)=\nn 
&=& \ell\uii t_R\ci \lsr\left(\ell^{(2^{\prime})} s_R\ci \ltr(a
\ell^{(3^{\prime})} )\ell\ui \right)\ot \ell^{(3)} \ot
S^{-1}\left(\ell^{(1^{\prime})}\right)=\nn  
%&=& t_R\ci \ltr (a \ell^{(3)^{\prime}} ) \ell\uii t_R\ci
%\lsr\left(\ell^{(2)^{\prime}} \ell\ui \right)\ot \ell^{(3)} \ot
%S^{-1}\left(\ell^{(1)^{\prime}}\right)=\nn 
&=& (\ltr\rd a \ell^{(3^{\prime})})\ell\uii t_R\ci
\lsr\left(\ell^{(2^{\prime})} \ell\ui \right)\ot \ell^{(3)} \ot
S^{-1}\left(\ell^{(1^{\prime})}\right)=\nn 
%&=& a\ui \ell^{(3)^{\prime}} s_R\ci \ltr (a\uii\ell^{(4)^{\prime}})
%\ell\uii t_R\ci
%\lsr\left(\ell^{(2)^{\prime}} \ell\ui \right)\ot \ell^{(3)} \ot
%S^{-1}\left(\ell^{(1)^{\prime}}\right)=\nn 
&=& a\ui\left( \lsr\ru \ell^{(2^{\prime})} s_R\ci \ltr
(a\uii\ell^{(3^{\prime})})  \ell\ui\right) \ot \ell^{(2)} \ot
S^{-1}\left(\ell^{(1^{\prime})}\right)=\nn 
%&=& a\ui s_R\ci \lsr \left( \ell^{(2)^{\prime}} s_R\ci \ltr
%(a\uii\ell^{(3)^{\prime}})  \ell\ui\right) \ot \ell^{(2)} \ot
%S^{-1}\left(\ell^{(1)^{\prime}}\right)=\nn 
&=& a\ui \ot \ell\uii t_R\ci \lsr \left( \ell^{(2^{\prime})} s_R\ci \ltr
(a\uii\ell^{(3^{\prime})})  \ell\ui\right) \ot
S^{-1}\left(\ell^{(1^{\prime})}\right)=\nn 
%&=& a\ui \ot \left(\lsr\lu[(\ltr\ld a\uii)\rd
%\ell^{(2)^{\prime}}]\right) \ru \ell  \ot
%S^{-1}\left(\ell^{(1)^{\prime}}\right)=\nn 
&=& a\ui\ot S^{-1}\left(S(a\uii)\uii\right) \ot
S^{-1}\left(S(a\uii)\ui\right). \hspace{1cm}\qed 
\nonumber\eea
One defines a {\em
non-degenerate right integral} in a left bialgebroid ${\cal
A}_L=(A,L,s_L,t_L,\gamma_L,\pi_L)$ as a non-degenerate left integral
in the right bialgebroid $({\cal A}_L)^{op}$. That is as
an element $\err\in A$ satisfying
\bea &i)& {both\  maps}\ {_L\err}\  {and}\  {\err_L}\ {are\
bijective}\nn
     &ii)& \err\di a\ot \err\dii=\err\di \ot \err \dii[\err\lu (a\ru
{\rho_*})]\nn
     && \err\di\ot \err\dii a = \err\di[\err\ld (a\rd {_* \rho})]\ot
\err\dii
\nonumber \eea
for all $a\in A$ and ${_*\rho}\colon = {_L \err}^{-1}(1_A)$ and
${\rho_*}\colon = {\err_L}^{-1}(1_A)$.
The application of Theorem  \ref{LS/2} to the right
bialgebroid $({\cal A}_L)^{op}$ implies that a left bialgebroid
${\cal A}_L$ 
possessing a non-degenerate right integral $\err$ can be made a Hopf
algebroid with the antipode $S(a)\colon = \err\lu(a\ru {\rho_*})$.

As an application of Theorem  \ref{LS/2} we sketch a different
derivation of the Hopf algebroid symmetry of an abstract depth 2
Frobenius extension -- obtained in \cite{FD2}. 

Let ${\cal C}$ be an additive 2-category closed under the direct sums and
subobjects of 1-morphisms. By an abstract extension we mean a
1-morphism $\i$ in $\cal C$ that possesses a left  dual $\ib$. This
means the existence of 2-morphisms 
\be ev_L\in \c^2({\bar{\iota}}\times {{\iota}},s_0(\i)) \qquad
    coev_L\in \c^2(t_0(\i),{{\iota}}\times {\bar{\iota}}) \lb{lint}\ee
satisfying the relations
\bea (\iota\times ev_L)\ci (coev_L\times \iota) &=&\iota \nn 
     (ev_L\times {\bar{\iota}})\ci (\ib \x coev_L) &=&\ib
\nonumber\eea
where $s_0(\i)$ and $t_0(\i)$ are the source and target 0-morphisms of
$\i$, respectively, and $\x$ stands for the horizontal product and
$\ci$ for the vertical product.
The 1-morphism $\i$  satisfies the {\em left depth
2} (or D2 for short) condition if $\i\x\ib \x\i$ is a direct summand in
a finite direct sum of copies of $\i$'s \cite{Sz}. This means the
existence of 
finite sets of 2-morphisms $\beta_i\in {\cal C}^2(\i\x\ib\x\i,\i)$ and
$\beta^{\prime}_i\in {\cal C}^2(\i,\i\x\ib\x\i)$ satisfying $\sum_i
\beta^{\prime}_i\ci \beta_i=\i\x\ib\x\i$. By the Theorem 3.5 in
\cite{Sz} in this case the ring of 2-morphisms: $A\colon = {\cal
C}^2(\i\x\ib,\i\x\ib)$ carries a left  bialgebroid structure ${\cal
A}_L$ over the base $L\colon ={\cal C}^2(\i,\i)$. The
structural maps of ${\cal A}_L$ are explicitly given in \cite{FD2}:
\bea s_L(l)&=&l\x \ib \nn
     t_L(l)&=&\i \x [({ev}_L\x\ib)\ci (\ib\x l \x \ib)\ci (\ib \x
{coev}_L)]\nn 
     \pi_L(a)&=&(\i \x {ev}_L)\ci (a \x \i)\ci ({coev}_L \x \i)\nn
     \gamma_L(a) &=& (\i\x {ev}_L\x \ib)\ci (a\x\i\x\ib)\ci (\i\x
{ev}_L\x \ib \x\i\x\ib)\ci (\i\x\ib \x \beta^{\prime}_i\x\ib)\ci
(\i\x\ib\x {coev}_L)\ot\nn
&&(\beta_i\x \ib)\ci (\i\x\ib\x {coev}_L)
\nonumber\eea
for $l\in L$ and $a\in A$. Summation on the index $i$ is understood.

The 1-morphism $\i$ satisfies the {\em right D2} condition if  $\ib\x\i
\x\ib$ is a direct summand in 
a finite direct sum of copies of $\ib$'s. This means the existence of
finite sets of 2-morphisms ${\hat \beta_i}\in {\cal
C}^2(\ib\x\i\x\ib,\ib)$ and 
${\hat \beta}^{\prime}_i\in {\cal C}^2(\ib,\ib\x\i\x\ib)$ satisfying $\sum_i
{\hat \beta}^{\prime}_i\ci {\hat \beta}_i=\ib\x\i\x\ib$. In this case
the ring $B\colon = {\cal 
C}^2(\ib\x\i,\ib\x\i)$ carries a right bialgebroid structure ${\cal
B}_R$ over the base $R\colon ={\cal C}^2(\i,\i)$. The
structural maps of ${\cal B}_R$ read as
\bea s_R(r)&=&\ib \x r \nn
     t_R(r)&=&[({ev}_L\x \ib)\ci (\ib\x r\x\ib)\ci (\ib\x {coev}_L)]\x
\i\nn
     \pi_R(b)&=&(\i \x {ev}_L)\ci (\i \x b)\ci ({coev}_L \x \i)\nn
     \gamma_R(b) &=& (\ib\x\i\x {ev}_L)\ci (\ib\x\i\x {\hat \beta}_i\x\i)\ci
(\ib\x {coev}_L \x \i \x \ib \x \i)\ci (b\x\ib\x\i)\ci (\ib\x
{coev}_L\x\i) \ot\nn
&&(\ib\x\i\x {ev}_L)\ci ({\hat \beta}^{\prime}_i\x\i)
\nonumber\eea
for $r\in R$ and $b\in B$. If $\iota$ satisfies both the left and the
right D2 conditions then the bialgebroids ${\cal A}_L$ and ${\cal B}_R$
are duals.

Let us assume that $\i$ is a {\em Frobenius} 1-morphism that is its dual
$\ib$ is two-sided. This means the existence of further 2-morphisms 
\be ev_R\in \c^2(\iota\times {\bar{\iota}},t_0(\i)) \qquad  coev_R\in
\c^2(s_0(\i),{\bar{\iota}}\times \iota) \lb{rint}\ee
satisfying the relations
\bea (ev_R\x \i)\ci (\i\x coev_R) &=&\i\nn
     (\ib \x ev_R)\ci (coev_R\x \ib) &=&\ib.
\nonumber\eea
Under this assumption the left- and  right D2 conditions become
equivalent (the 2-morphisms ${\hat \beta}_i$ and ${\hat
\beta}^{\prime}_i$ can be expressed in terms of $\beta_i$,
$\beta^{\prime}_i$ and the 2-morphisms (\ref{lint}-\ref{rint}).)
One checks that in the case of a D2 Frobenius 1-morphism $\i$ the element 
\be \err \colon = {coev}_L\ci {ev}_R \lb{re}\ee
of $A$ is a non-degenerate right integral. It leads to the antipode
\[S_A(a)=(\i\x\ib\x {ev}_R)\ci (\i\x\ib\x\i\x {ev}_L\x\ib)\ci 
(\i\x\ib\x a \x\i\x\ib)\ci (\i\x {coev}_R\x\ib\x\i\x\ib)\ci
({coev}_L\x\i\x\ib) \]
obtained by different methods in \cite{FD2}.
Also the element 
\be \ell\colon =  {coev}_R\ci {ev}_L \lb{ell}\ee
of $B$ is a non-degenerate left integral, leading to the antipode
\[ S_B(b)=(\ib\x\i\x{ev}_L)\ci (\ib\x\i\x \ib\x {ev}_R\x\i)\ci 
(\ib\x\i\x b \x\ib\x\i)\ci (\ib\x {coev}_L\x\i\x\ib\x\i)\ci
({coev}_R\x\ib\x\i). \]
As a matter of fact $S_A(\err)=\err$ and $S_B(\ell)=\ell$, that is
both $\err$ and $\ell$ turn out to be two-sided non-degenerate
integrals.  

We want to emphasize that the non-degenerate integrals (\ref{re}) and
(\ref{ell}) 
are non-unique. Another choice of the non-degenerate integrals leads
to other antipodes and other corresponding right and left bialgebroid
structures on $A$ and $B$, respectively. 

The name `abstract extension' is motivated by the most important
example. For an extension $N\to M$ of rings the forgetful functor
$\Phi$ of right modules ${\cal M}_M\to {\cal M}_N$ is a 1-morphism in
the 2-category of categories. It possesses a  left dual: the induction
functor. The 1-morphism $\Phi$  is left/right D2 
and Frobenius if and only if the 
extension $N\to M$ is left/right D2 and Frobenius, respectively. 
As it is shown in \cite{KSz} in the case of a D2 ring extension $N\to
M$ the above ring $A$ is isomorphic to the endomorphism ring ${\rm
End}({_NM_N})$ and $B$ is the center $(M\ot_N M)^N$. Now for a
Frobenius extension $N\to M$ let us {\em fix} a Frobenius system
that is an $N-N$ bimodule map $\psi:{_NM_N}\to N$ and its quasi-basis
$y_i\ot x_i \in M_N\ot {_N M}$. Then
the non-degenerate integrals (\ref{re}) and (\ref{ell}) are
$\err=\psi$ and $\ell=y_i\ot x_i$, respectively. The construction of
the corresponding antipodes gives the Hopf algebroid structure on $A$
that is discussed in the Section 3 of \cite{BKSz} on different grounds,
and its dual -- ${\cal A}_*^{\err}$ -- on $B$.

\end{document}